\documentclass[aap,seceqn,nameyear,MSNbibl,dvips]{arximspdf}

\doi{10.1214/10-AAP750}
\volume{21}
\issue{6}
\pubyear{2011}
\firstpage{2447}
\lastpage{2482}

\makeatletter

\newtheorem{theorem}{Theorem}
\newtheorem{lemma}{Lemma}[section]
\newtheorem{corollary}{Corollary}

\newcommand{\ep}{\varepsilon}

\newcommand{\bbr}{\mathbb{R}}
\newcommand{\bbz}{\mathbb{Z}}
\newcommand{\intt}{\int^t_0}
\newcommand{\eqp}{\stackrel{P}{\rightarrow}}
\newcommand{\ld}{\ldots}

\makeatother

\begin{document}
\begin{frontmatter}

\title{Asymptotic behavior of Aldous' gossip process\thanksref{T1}}
\runtitle{Aldous' gossip process}

\thankstext{T1}{Supported in part by NSF Grant DMS-07-04996 from the
probability program.}

\begin{aug}
\author[A]{\fnms{Shirshendu} \snm{Chatterjee}\ead[label=e1]{sc499@cornell.edu}} and
\author[B]{\fnms{Rick} \snm{Durrett}\corref{}\ead[label=e2]{rtd@math.duke.edu}}
\runauthor{S. Chatterjee and R. Durrett}
\affiliation{Cornell University and Duke University}
\address[A]{School of Operations Research\\
\quad and Information Engineering\\
Department of Mathematics\\
Cornell University\\
Ithaca, New York 14853\\
USA\\
\printead{e1}} 
\address[B]{Mathematics Department\\
Duke University\\
Box 90320\\
Durham, North Carolina 27708-0320\\
USA\\
\printead{e2}}
\end{aug}

\received{\smonth{5} \syear{2010}}
\revised{\smonth{9} \syear{2010}}

%
\begin{abstract}
Aldous [(2007) Preprint] defined a gossip process in
which space is a discrete $N \times N$ torus, and the state of the
process at time $t$ is the set of individuals who know the information.
Information spreads from a~site to its nearest neighbors at rate $1/4$
each and at rate $N^{-\alpha}$ to a site chosen at random from the
torus. We will be interested in the case in which $\alpha< 3$, where
the long range transmission significantly accelerates the time at which
everyone knows the information. We prove three results that precisely
describe the spread of information in a~slightly simplified model on
the real torus. The time until everyone knows the information is
asymptotically $T=(2-2\alpha/3) N^{\alpha/3} \log N$. If~$\rho_s$ is
the fraction of the population who know the information at time $s$ and
$\ep$ is small then, for large $N$, the time until $\rho_s$ reaches~%
$\ep$ is $T(\ep) \approx T + N^{\alpha/3} \log (3\ep/M)$, where $M$ is
a random variable determined by the early spread of the information.
The value of $\rho_s$ at time $s = T(1/3) + t N^{\alpha/3}$ is almost a
deterministic function $h(t)$ which satisfies an odd looking
integro-differential equation. The last result confirms a heuristic
calculation of Aldous.
\end{abstract}

%

\end{frontmatter}

\section{Introduction}\label{intro}

We study a model introduced by \citet{Ald07} for the spread of
gossip and other more economically useful information. His paper considers
various game theoretic aspects of random percolation of information
through networks. Here we concentrate on one small part, a
first passage percolation model with nearest neighbor and long-range jumps
introduced in his Section~6.2. The work presented here is also related to
work of \citet{FilMau04} and \citet{CanMarMon06},
who considered the impact of long-range dispersal on the spread of epidemics
and invading species.

Space is the discrete torus $\Lambda(N) = (\bbz\bmod N)^2$.
The state of the process at time $t$ is $\xi_t \subset\Lambda(N)$,
the set of individuals who know the information at time $t$.
Information spreads from $i$ to $j$ at rate
\[
\nu_{ij} = \cases{ 1/4, &\quad if $j$ is a (nearest) neighbor of $i$,\cr
\lambda_N/N^2, &\quad if not.}
\]
If $\lambda_N=0$, this is ordinary first passage percolation on the
torus. If we
start with $\xi_0 = \{(0,0)\}$, then the shape theorem for
nearest-neighbor\vadjust{\goodbreak} first passage percolation,
see \citet{CoxDur81} or \citet{Kes86}, implies that until the process
exits $(-N/2,N/2)^2$, the radius of the set $\xi_t$ grows linearly
and~$\xi_t$ has an asymptotic shape. From this we see that if $\lambda
_N=0$, then there is a~constant $c_0$ so that the time $T_N$, until everyone knows the
information, satisfies
\[
\frac{T_N}{N} \eqp c_0,
\]
where $\eqp$ denotes convergence in probability.

To simplify things, we will remove the randomness from the nearest
neighbor part of the process, and formulate it on the (real) torus
$\Gamma(N) = (\bbr\bmod N)^2$. One should be able to prove a similar
result for the first passage percolation model but there are two difficulties.
The first and easier to handle is that the limiting shape is not round.
The second
and more difficult issue is that the growth is not deterministic but
has fluctuations.
One should be able to handle both of these problems, but the proof is
already long enough.

We consider what we call the ``balloon process,'' in which
the state of the process at time $t$
is $\mathcal{C}_t \subset\Gamma(N)$. It starts with one ``center''
chosen uniformly from the torus at
time 0. When a center is born at $x$, a disk with radius 0 is put
there, and its radius grows deterministically as $r(s) = s/\sqrt{2\pi
}$, so that the
area of the disk at time $s$ after its birth is $s^2/2$. If the area
covered at time $t$ is $C_t$, then births of new centers occur at
rate $\lambda_NC_t$. The location of each new center is chosen
uniformly from the torus. If the new point lands at $x\in\mathcal{
C}_t$, it will never contribute anything to the growth of the set,
but we will count it in the total number of centers, which we denote by
$\tilde X_t$.

Before turning to the details of our analysis we would like to point
out that a related balloon process was used by Barbour and Reinert
(\citeyear{BR01}) in their study of distances on the small world graph.
Consider a
circle of radius $L$ and introduce a Poisson mean $\rho L/2$ number of
chords with length~0 connecting randomly chosen points on the circle.
To study the distance between a fixed point $O$ and a point chosen at
random one wants to examine $S(t) = \{ x \dvtx\operatorname{dist}(O,x)
\le t \}$.
If we ignore overlaps and let $M(t)$ be the number of intervals in
$S(t)$ then $S'(t) = 2M(t)$ and $M(t)$ is a Yule process with births at
rate $2\rho M(t)$ due to the interval ends encountering points in the
Poisson process of chords. This a balloon process in which the new
births come from the boundaries. As in our case one first studies the
growth of the ballon process and then estimates the difference from the
real process to prove the desired result. There are interesting
parallels and differences between the two proofs, see Section 5.2 of
Durrett (\citeyear{D07}) for a proof.

Here we will be concerned with $\lambda_N = N^{-\alpha}$. To begin
we will get rid of trivial cases. If the diameter of $\mathcal{C}_t$
grows linearly,
then $\int_0^{c_0 N} C_t \,dt = O(N^3)$.
So if $\alpha> 3$, with probability tending to 1 as $N$ goes to
$\infty$, there is no long range jump before
the initial disk covers the entire torus, and the time $T_N$ until the
entire torus is covered
satisfies
\[
\frac{T_N}{N} \eqp c_1 \qquad \mbox{where } c_1=\sqrt{\pi}.
\]
If $\alpha=3$, then with probabilities bounded
away from 0, (i) there is no long range jump and $T_N \approx
c_1N$, and (ii) there is one
that lands close enough to $(N/2,N/2)$ to make $T_N \le(1-\delta)
Nc_1$. Using
$\Rightarrow$ for weak convergence, this suggests that
\setcounter{theorem}{-1}
\begin{theorem}\label{theo0}
When $\alpha=3$, $T_N/N \Rightarrow$ a random limit concentrated on
$[0,c_1]$ and with an atom at
$c_1$.
\end{theorem}
\begin{pf}
Suppose without loss of generality that the initial center is at~0, and
view the torus as
$(-N/2,N/2]^2$. The key observation is that the set-valued process $\{
\mathcal{C}_{Nt}/N , t\ge0 \}$
converges to a limit $\mathcal{D}_t$. Before the first long-range
dispersal, the state of $\mathcal{D}_t$ is the intersection of
the disk of radius $t/\sqrt{2\pi}$ with $(-1/2,1/2]^2$. Long range
births occur at rate equal to the
area of $\mathcal{D}_t$ and are dispersed uniformly. Since the
distance from $(0,0)$ to $(1/2,1/2)$ is $1/\sqrt{2}$,
if there are no long range births before time $c_1=\sqrt{\pi}$ or if
all long range births land inside $\mathcal{D}_t$
then the torus is covered at time $c_1$. Computing the distribution of
the cover time when it is $< c_1$ is
complicated, but the answer is a continuous functional of the limit
process, and standard weak convergence results
give the result.
\end{pf}

For the remainder of the paper we suppose $\lambda_N = N^{-\alpha}$
with $\alpha<3$. The overlaps between disks in $\mathcal{C}_t$ pose a
difficulty in analyzing the process, so we begin by studying a
simpler ``balloon branching process'' $\mathcal{A}_t$, in which $A_t$ is
the sum of the areas of all of the disks at time $t$, births of new
centers occur at rate $\lambda_NA_t$, and the location of each new
center is chosen uniformly from the torus. Let $X_t$ be the number
of centers at time $t$ in $\mathcal{A}_t$.

Suppose we start $\mathcal{C}_0$ and $\mathcal{A}_0$ from the same
randomly chosen point.
The areas $C_t=A_t$ until the time of the first birth, which can be
made to be the same in the two processes.
If we couple the location of the new centers at that time, and continue
in the obvious way letting $\mathcal{C}_t$
and $\mathcal{A}_t$ give birth at the same time with the maximum rate
possible, to the same place
when they give birth simultaneously, and letting\vadjust{\goodbreak}
$\mathcal{A}_t$ give birth by itself otherwise, then we will have
%
%
\begin{equation}\label{couple}
\mathcal{C}_t \subset\mathcal{A}_t,\qquad C_t \le A_t,\qquad
\tilde X_t \le X_t \qquad\mbox{for all $t\ge0$.}
\end{equation}

$X_t$ is a Crump--Mode--Jagers branching process, but saying these
words does not magically solve our problems. Define the length
process $L_t$ to be $\sqrt{2\pi}$ times the sum of the radii of all
the disks at time $t$.
%
%
\begin{eqnarray} \label{LA}
L_t &=&\int_0^t (t-s) \,dX_s = \int_0^t X_s \,ds, \nonumber\\[-8pt]\\[-8pt]
A_t &=&\int_0^t \frac{(t-s)^2}{2} \,dX_s = \int_0^t L_s \,ds.
\nonumber
\end{eqnarray}
Here and later we use $\intt$ for integration over the closed
interval $[0,t]$, that is, we include the contribution from the atom in
$dX_s$ at 0 ($X_0=1$ while $X_s=0$ for $s<0$). For the second
equality on each line integrate by parts or note that $L_t'=X_t$ and
$A_t'=L_t$. Since $X_t$ increases by 1 at rate $\lambda_NA_t$,
$(X_t,L_t,A_t)$ is a Markov process.

To simplify formulas, we will often drop the subscript $N$ from
$\lambda_N$. For comparison with $C_t$, the parameter $\lambda$ is
important, but in the analysis of $A_t$ it is not. If we let
%
%
\begin{equation}\label{scale}\qquad
X^1_t = X(t\lambda^{-1/3}),\qquad
L^1_t = \lambda^{1/3} L(t\lambda^{-1/3}),\qquad
A^1_t = \lambda^{2/3} A(t\lambda^{-1/3}),
\end{equation}
then $(X_t^1,L^1_t,A^1_t)$ is the process with $\lambda=1$.

To study the growth of $A_t$, first we will compute the means of
$X_t$, $L_t$ and~$A_t$. Let $F(t) = \lambda t^3/3!$. Using the
independent and identical behavior of all the disks in $\mathcal{A}_t$
it is easy to show that (see the proof of Lemma \ref{XLAlem})
\[
EX_t = 1 + \int_0^t EX_{t-s} \,dF(s).
\]
Solving the above renewal equation and using (\ref{LA}), we can show
%
%
\begin{eqnarray}\label{mean}
EX_t &=& \sum_{k=0}^\infty F^{*k}(t) = V(t) = \sum_{k=0}^\infty\frac
{\lambda^k t^{3k}}{(3k)!}, \nonumber\\
EL_t &=& \sum_{k=0}^\infty\frac{\lambda^k
t^{3k+1}}{(3k+1)!},\\
EA_t &=& \sum_{k=0}^\infty
\frac{\lambda^k t^{3k+2}}{(3k+2)!}. \nonumber
\end{eqnarray}
To evaluate $V(t)$ we note that $V'''(t)\!=\!\lambda V(t)$ with
$V(0)\!=\!1, V'(0)\!=\!V''(0)\!=\!0$, so
%
%
\begin{equation}\label{Vtdef}
V(t)=\tfrac{1}{3}[\exp(\lambda^{1/3}t)
+\exp(\lambda^{1/3}\omega t)
+\exp(\lambda^{1/3}\omega^2 t)].
\end{equation}
Here $\omega=(-1+i\sqrt{3})/2$ is one
of the complex cube roots of 1 and
$\omega^2=(-1-i\sqrt{3})/2$ is the other. Note that each
of $\omega$ and $\omega^2$ has real part $-1/2$. So the second and
third terms in (\ref{Vtdef}) go to 0 exponentially fast.

If $\mathcal{F}_s=\sigma\{X_r, L_r, A_r\dvtx r\le s\}$, then
%
%
\begin{equation}\label{infgen}
\frac{d}{dt} E\left.\left.\left[
\matrix{X_t \cr L_t \cr A_t}
\right| \mathcal{F}_s \right] \right|_{t=s} =
\pmatrix{
0 & 0 & \lambda\cr
1 & 0 & 0 \cr
0 & 1 & 0
}
\left[\matrix{X_s \cr L_s \cr A_s}\right].
\end{equation}
Let $Q$\vspace*{1pt} be the matrix in (\ref{infgen}). By
computing the determinant of $Q-\eta I$ it is easy to see that $Q$
has eigenvalues $\eta= \lambda^{1/3}, \omega\lambda^{1/3},
\omega^2 \lambda^{1/3}$, and $e^{-\eta t} ( X_t + \eta L_t + \eta^2
A_t )$ is a (complex) martingale. To treat the three martingales
separately, let
\begin{eqnarray*}
I_t &=& X_t + \lambda^{1/3}L_t + \lambda^{2/3}A_t, \qquad
M_t=\exp(-\lambda^{1/3}t) I_t,\\
J_t &=& X_t + (\omega\lambda^{1/3})L_t + (\omega\lambda
^{1/3})^2A_t, \qquad
\tilde J_t=\exp(-\omega\lambda^{1/3}t) J_t,\\
K_t &=& X_t + (\omega^2\lambda^{1/3})L_t + (\omega^2\lambda
^{1/3})^2A_t, \qquad
\tilde K_t=\exp(-\omega^2\lambda^{1/3}t) K_t,
\end{eqnarray*}
so that $M_t$ is the real martingale, and $\tilde J_t$ and $\tilde
K_t $ are the complex ones.
\begin{theorem} \label{th1}
$\{M_t\dvtx t\ge0\}$ is a positive square integrable martingale with
respect to the filtration $\{ \mathcal F_t\dvtx t\ge0\}$. $EM_t=M_0=1$.
\begin{eqnarray*}
& EM_t^2 =\frac{8}{7}-\frac{1}{3} \exp(-\lambda^{1/3}t)
+ O\bigl(\exp(-5\lambda^{1/3}t/2)\bigr),& \\
& E |\tilde J_t|^2,\ E|\tilde K_t|^2 = \frac{1}{6} \exp(2\lambda
^{1/3}t) +
O\bigl( \exp(\lambda^{1/3}t/2)\bigr).&
\end{eqnarray*}
If we let $M = \lim_{t\to\infty} M_t$, then $P(M>0)=1$ and
\[
\exp(-\lambda^{1/3}t)X_t,\mbox{ }
\lambda^{1/3}\exp(-\lambda^{1/3}t)L_t,\mbox{ }
\lambda^{2/3}\exp(-\lambda^{1/3}t)A_t \to M/3
\]
a.s. and in $L^2$. The distribution of $M$ does not depend on $\lambda$.
\end{theorem}

The last result follows from (\ref{scale}), which with (\ref{LA})
explains why the three quantities
converge to the same limit.
The key to the proof of the convergence results is to note that
$1+\omega+\omega^2=0$ implies
\begin{eqnarray*}
3X_t &=& I_t + J_t + K_t, \\
3\lambda^{1/3} L_t &=& I_t + \omega^2 J_t + \omega K_t, \\
3\lambda^{2/3} A_t &=& I_t + \omega J_t + \omega^2 K_t. \nonumber
\end{eqnarray*}
The real parts of $\omega$ and $\omega^2$ are $-1/2$.
Although the results for $E|\tilde J_t|^2$ and $E|\tilde K_t|^2$ show
that the martingales
$\tilde J_t$ and $\tilde K_t$ are not $L^2$ bounded, it is easy to
show that $\exp(-\lambda^{1/3}t)J_t$ and
$\exp(-\lambda^{1/3}t) K_t \to0$ a.s. and in $L^2$, and
Theorem \ref{th1} then follows from $M_t = \exp(-\lambda
^{1/3}t)
I_t \to M$.

Recall that $\lambda_N = N^{-\alpha}$ and let
%
%
\begin{eqnarray}\label{a}
a(t) &=& (1/3) N^{2\alpha/3} \exp( N^{-\alpha/3} t),\qquad
l(t)=N^{-\alpha/3}a(t),\nonumber\\[-8pt]\\[-8pt]
x(t) &=& N^{-2\alpha/3} a(t),\nonumber
\end{eqnarray}
so that
$A_t/a(t), L_t/l(t), X_t/x(t) \to M$ a.s. Let
%
%
\begin{equation}\label{S}
S(\ep) = N^{\alpha/3}[
(2-2\alpha/3) \log N + \log(3\ep) ],
\end{equation}
so $a(S(\ep))=\ep N^2$. Let
%
%
\begin{equation} \label{sigtau}
\sigma(\ep) = \inf\{ t \dvtx A_t \ge
\ep N^2 \} \quad\mbox{and}\quad \tau(\ep) = \inf\{ t \dvtx C_t \ge\ep
N^2 \}.
\end{equation}
The first of these is easy to study.
\begin{theorem} \label{th2}
If $0<\ep< 1$, then as $N\to\infty$
\[
N^{-\alpha/3} \bigl(\sigma(\ep) - S(\ep)\bigr) \eqp- \log(M).
\]
The coupling in (\ref{couple}) implies $\tau(\ep)\ge\sigma(\ep)$. In
the other direction, for any $\gamma>0$
\[
\limsup_{N\to\infty} P\bigl[ \tau(\ep) > \sigma\bigl((1+\gamma)\ep\bigr) \bigr] \le
P\bigl( M \le(1+\gamma)\ep^{1/3} \bigr) + 11\frac{\ep^{1/3}}
{\gamma}.
\]
\end{theorem}

The last result implies that for $\ep<1$
%
%
\begin{equation} \label{tauLLN}
\tau(\ep) \sim
(2-2\alpha/3) N^{\alpha/3}\log N.
\end{equation}
Our next goal
is to obtain more precise information about $\tau(\ep)$ and about
how $|C_t|/N^2$ increases from a small positive level to reach 1.

The first result in Theorem \ref{th2} shows that
$(\sigma(\ep)-S(\ep))/N^{\alpha/3}$ is determined by the random
variable $M$ from Theorem \ref{th1}, which in turn is determined by
what happens early in the growth of the branching balloon process.
Let
%
%
\begin{equation}
\label{R}
R = N^{\alpha/3}[(2-2\alpha/3)\log N - \log(M)],
\end{equation}
$R$ is defined so that $a(R) = (1/3) N^2 /M$, and hence
$A_R/N^2 \eqp1/3$. Define
%
%
\begin{equation}\label{psiWI}\qquad
\psi(t)\equiv
R+N^{\alpha/3}t,\qquad W\equiv\psi(\log(3\ep))\quad \mbox{and}\quad
I_{\ep,t}=[\log(3\ep) , t]
\end{equation}
for $\log(3\ep) \le t$. $W$ is
defined so that $a(W)=\ep N^2/M$ and hence $A_W/N^2 \eqp\ep$. The
arguments that led to Theorem \ref{th2} will show that if $\ep$ is
small then~$C_W/A_W$ is close to 1 with high probability.

To get a lower bound on the growth of $C_t$ after time $W$ we
declare that the centers in $\mathcal{C}_W$ and $\mathcal{A}_W$ to be
generation 0 in $\mathcal{C}_t$ and $\mathcal{A}_t$, respectively, and we
number the succeeding generations in the obvious way, a center born
from an area of generation $k$ is in generation $k+1$. For $t\ge
\log(3\ep)$, let $C_{W,\psi(t)}^k$ and $A_{W,\psi(t)}^k$ denote\vadjust{\goodbreak} the
areas covered at time $\psi(t)$ by respective centers of generations
$j\in\{0, 1, \ldots, k\}$ and let
%
%
\begin{eqnarray}\label{gfdef}
g_{0}(t)&=&\ep\biggl[1+\bigl(t-\log(3\ep)\bigr)+\frac{(t-\log(3\ep
))^2}{2}\biggr],\nonumber\\[-8pt]\\[-8pt]
f_0(t)&=&g_0(t)-\ep^{7/6}.\nonumber
\end{eqnarray}
To explain these
definitions, we note that Lemma \ref{B0bounds} will show that for
any~$t$, there is an $\ep_0=\ep_0(t)$ so that for any $0 < \ep<
\ep_0$
\begin{eqnarray*}
\lim_{N\to\infty} P\Bigl(\sup_{s\in I_{\ep,t}}
\bigl|N^{-2}A^0_{W,\psi(s)}-g_0(s)\bigr|>\eta\Bigr)&=&0\qquad
\mbox{for any $\eta>0$},\\
P\Bigl( \inf_{s\in I_{\ep,t}} N^{-2} \bigl(C^0_{W,\psi(s)} -
A^0_{W,\psi(s)}\bigr) < - \ep^{7/6} \Bigr) &\le& P( M < \ep^{1/3} ) +
\ep^{1/12}.
\end{eqnarray*}
Since $C_{W,\psi(t)}^0 \le A_{W,\psi(t)}^0$, if $\ep$ is small, with
high probability $g_0(t)$ and $f_0(t)$ provide upper and lower
bounds, respectively, for $C_{W,\psi(t)}^0$.

To begin to improve these bounds we let
\[
f_1(t) = 1-\bigl(1-f_0(t)\bigr)\exp\biggl(-\int_{\log(3\ep)}^t
\frac{(t-s)^2}{2} f_{0}(s) \,ds\biggr),
\]
and define $g_1$ similarly. To explain this equation note that an $x
\notin C_{W,\psi(t)}^0$
will not be in $C_{W,\psi(t)}^1$ if and only if no generation 1 center
is born in the space--time cone
\[
K_{x,t}^\ep\equiv\bigl\{(y,s)\in\Gamma(N)\times[W,\psi(t)]\dvtx
|y-x| \le\bigl(\psi(t)-s\bigr)/\sqrt{2\pi}\bigr\}.
\]
Lemma \ref{f1lb} shows that for $0< \ep<\ep_0$ and $\delta>0$,
\[
\limsup_{N\to\infty} P\Bigl( \inf_{s\in I_{\ep,t}}
N^{-2}C^1_{W,\psi(s)} - f_1(s) < - \delta\Bigr) \le P( M <
\ep^{1/3} ) + \ep^{1/12}.
\]
To iterate this we will let
\[
f_{k+1}(t) = 1 - \bigl(1-f_{k}(t)\bigr) \exp\biggl(-\int_{\log(3\ep)}^t
\frac{(t-s)^2}{2}\bigl(f_k(s)-f_{k-1}(s)\bigr) \,ds\biggr)
\]
for $k\ge1$.
The difference $f_k(s)-f_{k-1}(s)$ in the integral comes from the
fact that a new point in generation $k+1$ must come from a point
that is in generation $k$ but not in generation $k-1$. Combining
these equations we have\looseness=-1
\begin{eqnarray*}
&&
1-f_{k+1}(t)\\
&&\qquad = \bigl(1-f_k(t)\bigr) \exp\biggl(-\int_{\log(3\ep)}^t \frac{(t-s)^2}{2}
\bigl(f_k(s)-f_{k-1}(s)\bigr) \,ds\biggr)\\
&&\qquad = \bigl(1-f_{k-1}(t)\bigr)\exp\Biggl(-\int_{\log(3\ep)}^t \frac{(t-s)^2}{2}
\sum_{l=k-1}^k\bigl(f_l(s)-f_{l-1}(s)\bigr) \,ds\Biggr)\cdots\\
&&\qquad = \bigl(1-f_0(t)\bigr)\exp\Biggl(-\int_{\log(3\ep)}^t
\frac{(t-s)^2}{2} \sum_{l=1}^k \bigl(f_l(s)-f_{l-1}(s)\bigr) + f_0(s)
\,ds\Biggr)
\end{eqnarray*}\looseness=0
so that
%
%
\begin{equation} \label{fkinteq}
f_{k+1}(t) =
1-\bigl(1-f_0(t)\bigr)\exp\biggl(-\int_{\log(3\ep)}^t \frac{(t-s)^2}{2}
f_{k}(s) \,ds\biggr).
\end{equation}
Since $f_1(t) \ge
f_0(t)$, letting $k\to\infty$, $f_k(t)\uparrow f_\ep(t)$, where
$f_\ep$ is the unique solution of
%
%
\begin{equation} \label{fepinteq}
f_\ep(t)=1-\bigl(1-f_0(t)\bigr)\exp\biggl(-\int_{\log(3\ep)}^t
\frac{(t-s)^2}{2} f_\ep(s) \,ds\biggr)
\end{equation}
with
$f_\ep(\log(3\ep))=\ep-\ep^{7/6}$. $g_k(t)$ and $g_\ep(t)$ are
defined similarly.

$g_\ep(t)$ and $f_\ep(t)$ provide upper and lower bounds on the
growth of $C_{\psi(t)}$ for $t \ge\log(3\ep)$.
To close the gap between these bounds we let $\ep\to0$.
\begin{lemma}\label{h}
For any $t<\infty$, if $I_{\ep,t}=[\log(3\ep),t]$, then as $\ep\to
0$,
\[
\sup_{s\in I_{\ep,t}} |f_\ep(s)-h(s)|\mbox{, }
\sup_{s\in I_{\ep,t}} |g_\ep(s)-h(s)| \to0
\]
for some nondecreasing $h$ with \textup{(a)} $\lim_{t\to-\infty} h(t) = 0$,
\textup{(b)} $\lim_{t\to\infty} h(t) = 1$,
{\renewcommand{\theequation}{c}
\begin{equation}
h(t) = 1-\exp\biggl(-\int_{-\infty}^t \frac{(t-s)^2}{2} h(s)
\,ds\biggr)
\end{equation}
}

\vspace*{6pt}

\vspace*{-\baselineskip}

\noindent and \textup{(d)} $0 < h(t) < 1$ for all $t$.
\end{lemma}

If one removes the 2 from inside the exponential, this is equation (36)
in \citet{Ald07}.
Since there is no initial condition, the solution is only unique up to
time translation.
\begin{theorem}\label{th3}
Let $h$ be the function in Lemma \ref{h}. For any $t<\infty$ and
$\delta>0$,
\[
\lim_{N\to\infty} P\Bigl(\sup_{s\le t} \bigl|N^{-2}C_{\psi(s)}-h(s)\bigr|
\le
\delta\Bigr)=1.
\]
\end{theorem}

This result shows that the displacement of $\tau(\ep)$ from
$(2-2\alpha/3)N^{\alpha/3} \log N$ on the
scale $N^{\alpha/3}$ is dictated by the random variable $M$ that
gives the rate of growth of the branching balloon process, and that once
$C_t$ reaches $\ep N^2$, the growth is deterministic.

The solution $h(t)$ never reaches 1, so we need a little more work to
show that
\begin{theorem} \label{th4} Let $T_N$ be the first time the torus is
covered. As $N\to\infty$
\[
T_N / (N^{\alpha/3} \log N) \eqp2-2\alpha/3.\vadjust{\goodbreak}
\]
\end{theorem}

The remainder of the paper is organized as follows. In Section \ref{sec2},
we prove the properties of $\mathcal{A}_t$ presented in Theorem
\ref{th1}. In Section \ref{sec3}, we prove the properties of the hitting
times s $\sigma(\ep)$ and $\tau(\ep)$ stated in Theorem \ref{th2}.
In Section \ref{sec4}, we prove the limiting behavior of $\mathcal{ C}_t$
mentioned in Theorem \ref{th3}. Finally in Section \ref{sec5}, we
prove Theorem \ref{th4}.

\section{Properties of the balloon branching process
$\mathcal{A}_t$}\label{sec2}

\begin{lemma}\label{conv}
$\intt s^m(t-s)^n \,ds=\frac{m!n!}{(m+n+1)!}t^{m+n+1}$.
\end{lemma}
\begin{pf} If you can remember the definition of the beta distribution,
this is trivial.
If you cannot then integrate by parts and use induction.
\end{pf}

Let $F(t)=\lambda t^3/3!$ for $t \ge0$, and $F(t)=0$ for $t<0$. Let
$V(t)=\sum_{k=0}^\infty F^{*k}(t)$, where $*k$ indicates the
$k$-fold convolution.
\begin{lemma}\label{V} If $\omega=(-1+i\sqrt{3})/2$, then
\[
V(t)=\sum_{k=0}^\infty\frac{\lambda^k t^{3k}}{(3k)!}
=\frac{1}{3}[\exp(\lambda^{1/3}t)+\exp
(\lambda^{1/3}\omega
t)+\exp(\lambda^{1/3}\omega^2 t)].
\]
\end{lemma}
\begin{pf}
We first use induction to show that
%
%
\begin{equation}\label{Fconv}
F^{*k}(t)= \cases{ \lambda^kt^{3k}/(3k)!, &\quad$t\ge0$,\cr
0, &\quad$t<0$.}
\end{equation}
This holds for $k=0, 1$ by our assumption. If the equality holds for
$k=n$, then using Lemma \ref{conv}
we have for $t \ge0$
\[
F^{*(n+1)}(t)=\int_0^t F^{*n}(t-s) \,dF(s) =\int_0^t
\frac{\lambda^n(t-s)^{3n}}{(3n)!}\frac{\lambda s^2}{2} \,ds
=\frac{\lambda^{n+1}t^{3n+3}}{(3n+3)!}.
\]
It follows by induction that $V(t)=\sum_{k=0}^\infty\lambda
^kt^{3k}/(3k)!$. To evaluate the sum
we note that setting $\lambda=1$, $U(t)=\sum_{k=0}^\infty
t^{3k}/(3k)!$ solves
\[
U'''(t)=U(t) \qquad\mbox{with $U(0)=1$ and $U'(0)=U''(0)=0$.}
\]
This differential equation has solutions of the from $e^{\gamma t}$,
where $\gamma^3=1$, that is, $\gamma=1, \omega$ and $\omega^2$. This
leads to the general solution
\[
U(t)=Ae^t+Be^{\omega t}+Ce^{\omega^2 t}
\]
for some constants $A, B, C$. Using the initial conditions for
$U(t)$ we have
\[
A+B+C=1,\qquad A+B\omega+C\omega^2=0,\qquad A+B\omega^2+C\omega=0.
\]
Since $1+ \omega+ \omega^2=0$, we have $A=B=C=1/3$. Since $V(t) =
U(\lambda^{1/3}t)$,
we have proved the desired result.
\end{pf}

Our next step is to compute the first two moments of $X_t, L_t$ and
$A_t$. For that we need the following lemma in addition to the
previous one.
\begin{lemma}\label{renewaleq}
Let $\{N_t\dvtx t\ge0\}$ be a Poisson process on $[0,\infty)$ with
intensity $\lambda(\cdot)$
and let $\Pi_t$ be the set of points at time $t$.
If $\{Y_t,Z_t\dvtx t\ge0\}$ are two complex valued stochastic processes
satisfying
\[
Y_t=y(t)+\sum_{s_i\in\Pi_t} Y^i_{t-s_i},\qquad Z_t=z(t)+\sum
_{s_i\in\Pi_t} Z^i_{t-s_i},
\]
where $(Y^i, Z^i)$, $i=1, 2, \ldots,$ are i.i.d. copies of $(Y,Z)$,
and independent of~$N$, then
\begin{eqnarray*}
EY_t &=&y(t)+\intt EY_{t-s}\lambda(s) \,ds, \\
E(Y_tZ_t) &=&(EY_t)(EZ_t)+\intt E(Y_{t-s}Z_{t-s})\lambda(s) \,ds.
\end{eqnarray*}
\end{lemma}
\begin{pf}
$N_t$ has Poisson distribution with mean $\Lambda_t=\intt\lambda(s) \,ds$.
Given $N_t=n$, the conditional distribution of $\Pi_t$ is same as the
distribution of $\{t_1, \ld, t_n\}$, where $t_1, \ldots, t_n$ are
i.i.d. from $[0,t]$ with
density $\beta(\cdot)=\lambda(\cdot)/\Lambda_t$. Hence
\[
E(Y_t|N_t)=y(t)+\sum_{i=1}^{N_t} EY^i_{t-t_i} =y(t)+N_t\int_0^t
EY_{t-s} \beta(s) \,ds,
\]
and taking expected values $EY_t=y(t)+\intt EY_{t-s}\lambda(s) \,ds$.

Similarly $EZ_t=z(t)+\intt EZ_{t-s}\lambda(s) \,ds$.
Using the conditional distribution of $\Pi_t$ given $N_t$,
\begin{eqnarray*}
E(Y_tZ_t|N_t)&=&y(t)z(t) + y(t) E\sum_{i=1}^{N_t} Z^i_{t-t_{i}} +z(t) E\sum
_{i=1}^{N_t} Y^i_{t-t_{i}}\\[-2pt]
&&{}+E\Biggl[\sum_{i=1}^{N_t} Y^i_{t-t_i}Z^i_{t-t_i}+\sum_{i\ne
j}Y^i_{t-t_i}Z^j_{t-t_j}\Biggr]\\[-2pt]
&=&y(t)z(t)+y(t)N_t\int_0^t EZ_{t-s} \beta(s) \,ds\\[-2pt]
&&{}  + z(t)N_t\int
_0^t EY_{t-s} \beta(s) \,ds+N_t\int_0^t E(Y_{t-s}Z_{t-s}) \beta(s) \,ds\\[-2pt]
&&{}+N_t(N_t-1)\int_0^t EY_{t-s} \beta(s) \,ds \int_0^t EZ_{t-s}\beta
(s) \,ds.
\end{eqnarray*}
Taking expectation on both sides and using
$EN_t(N_t-1)=\Lambda^2_t$, we get
\[
E(Y_tZ_t) = (EY_t)(EZ_t)+\intt E(Y_{t-s}Z_{t-s})\lambda(s) \,ds,
\]
which completes the proof.\vadjust{\goodbreak}
\end{pf}

Now we use Lemmas \ref{V} and \ref{renewaleq} to have the first
moments.
\begin{lemma} $E(X_t, L_t,A_t) = (V(t),V''(t)/\lambda, V'(t)/\lambda)$.
\label{XLAlem}
\end{lemma}
\begin{pf} Recall that $F(t)=\lambda t^3/3!$. In the balloon branching
process, the initial center $x$ gives
birth to new centers at rate $F'(t) = \lambda t^2/2$, and all the
centers behave independently and with the same distribution as the
one at~$x$. So
%
%
\begin{equation} \label{Xbreakup}
X_t=1+\sum_{s_i \in\Pi_t} X^i_{t-s_i},
\end{equation}
where
$\Pi_t \subset[0,t]$ is the set of times when new centers are born
in $\mathcal{A}_t$ and~$X^i$, $i=1, 2, \ldots,$ are i.i.d. copies of $X$,
and using Lemma \ref{renewaleq},
\[
EX_t =1 + \int_0^t EX_{t-s} \,dF(s).
\]
Using (4.5) from Chapter 3 of \citet{Dur10} and then (\ref{LA}):
%
%
\begin{eqnarray}\label{meanXLA}
EX_t &=& V(t) = \sum_{k=0}^\infty\frac{\lambda^k t^{3k}}{(3k)!},
\nonumber\\
EL_t &=&\int_0^t EX_s \,ds
=\sum_{k=0}^\infty\frac{\lambda^k t^{3k+1}}{(3k+1)!}, \\
EA_t &=& \int_0^t EL_s \,ds =\sum_{k=0}^\infty\frac{\lambda^k
t^{3k+2}}{(3k+2)!}. \nonumber
\end{eqnarray}
Since $V(t) = 1 + \sum_{k=0}^\infty\lambda^{k+1} t^{3k+3}/(3k+3)!$,
it is easy to see that $EA_t=V'(t)/\lambda$ and
$EL_t=V''(t)/\lambda$.
\end{pf}
\begin{lemma}\label{mart}
If $M_t
=\exp(-\lambda^{1/3}t)[X_t+\lambda^{1/3}L_t+\lambda^{2/3}A_t]$, then
$\{M_t\dvtx t\ge0\}$ is a square integrable martingale with respect to
the filtration $\{ \mathcal F_t\dvtx t\ge0\}$. $EM_t=1$ and
\[
EM_t^2=\tfrac{8}{7}-\tfrac{1}{3} \exp(-\lambda^{1/3}t) +
\theta_t \qquad\mbox{where } |\theta_t| \le\tfrac{4}{15}
\exp(-5\lambda^{1/3}t/2)
\]
and hence $(8/7) - EM_t^2 \le\exp(-\lambda^{1/3}t)$.
\end{lemma}
\begin{pf} Let $h(t,x,\ell,a) = \exp(-\lambda^{1/3}t)[x+\lambda
^{1/3}\ell+\lambda^{2/3}a]$, and let
$\mathcal{L}$ be the generator of the Markov process $(t,X_t,L_t,A_t)$.
Equation (\ref{infgen}) implies $\mathcal{L}h=0$, so $M_t$ is a martingale from
Dynkin's formula. $EM_t=EM_0=1$.

To compute $EM_t^2$ we use Lemma \ref{renewaleq} as follows. Let
$Y_t=Z_t=X_t+\lambda^{1/3}L_t+\lambda^{2/3}A_t$ and $g(t)\equiv
(EY_t)^2$. Since $EM_t=1$, $g(t)=\exp(2\lambda^{1/3}t)$.
Combining (\ref{LA})\vadjust{\goodbreak} and (\ref{Xbreakup}), letting $L_t^i = \int_0^t
X_s^i \,ds$ and $A_t^i = \int_0^t L_s^i \,ds, i=1, 2, \ldots,$ and
changing the variables $u=s-s_i$, we see that
\[
L_t=\int_0^t \biggl[1
+ \sum_{s_i \in\Pi_s} X_{s-s_i}^i
\biggr] \,ds = t + \sum_{s_i \in\Pi_t} \int_0^{t-s_i}X_u^i \,du =
t + \sum_{s_i \in\Pi_t}
L_{t-s_i}^i
\]
and hence
\[
A_t=\int_0^t \biggl[t
+ \sum_{s_i \in\Pi_s} L_{s-s_i}^i
\biggr] \,ds = t^2/2 + \sum_{s_i \in\Pi_t} \int_0^{t-s_i}L_u^i
\,du = t^2/2 + \sum_{s_i \in\Pi_t}
A_{t-s_i}^i.
\]
Thus all of $X_t, L_t$ and $A_t$ satisfy the hypothesis of Lemma
\ref{renewaleq} and so do~$Y_t$ and~$Z_t$, as they are linear
combinations of $X_t, L_t$ and $A_t$.
So applying Lem\-ma~\ref{renewaleq}
\[
EY_t^2=g(t)+\intt EY_{t-s}^2 \,dF(s).
\]
Solving the renewal equation using (4.8) in Chapter 3 of \citet{Dur10},
\[
EY_t^2=g*V(t)=\exp(2\lambda^{1/3}t)+\intt
\exp\bigl(2\lambda^{1/3}(t-s)\bigr) V'(s) \,ds,
\]
where $V=\sum_{k=0}^\infty F^{*k}$. To evaluate the integral we use
Lemma \ref{V} to conclude
\begin{eqnarray*}
&&\int_0^t \exp(-2\lambda^{1/3}s) V'(s) \,ds \\
&&\qquad=\frac{1}{3} \int_0^t \exp(-2\lambda^{1/3}s) \\
&&\qquad\quad\hspace*{21pt}{}\times\lambda^{1/3}[\exp(\lambda^{1/3}s)+\omega
\exp(\lambda^{1/3}\omega s)
+\omega^2 \exp(\lambda^{1/3}\omega^2 s)] \,ds \\
&&\qquad=\frac{1}{3}\biggl[ \frac{1}{1-2}\{
\exp(-\lambda^{1/3}t)-1\}
+\frac{\omega}{\omega-2}\bigl\{\exp\bigl((\omega-2)\lambda
^{1/3} t\bigr)-1\bigr\}\\
&&\qquad\quad\hspace*{112.1pt}{}
+\frac{\omega^2}{\omega^2-2}\bigl\{\exp\bigl((\omega
^2-2)\lambda^{1/3}
t\bigr)-1\bigr\}\biggr].
\end{eqnarray*}
Now using $1= -\omega-\omega^2$ and
$\omega^3=1$,
\[
1-\frac{\omega}{\omega-2}-\frac{\omega^2}{\omega^2-2}
=1-\frac{\omega^3 - 2\omega+\omega^3
-2\omega^2}{\omega^3-2\omega^2-2\omega^2+4}
=1-\frac{4}{7}=\frac{3}{7}.
\]
Since $\omega= (-1+i\sqrt{3})/2$ and $\omega^2 =
(-1-i\sqrt{3})/2$, the remaining error satisfies
\begin{eqnarray*}
3|\theta_t| &=&
\biggl|\frac{\omega}{\omega-2}\exp\bigl((\omega-2)\lambda^{1/3}
t\bigr)\biggr|
+ \biggl| \frac{\omega^2}{\omega^2-2}\exp\bigl((\omega
^2-2)\lambda^{1/3} t\bigr) \biggr| \\
&=& \biggl( \frac{1}{|\omega-2|} + \frac{1}{|\omega^2-2|} \biggr)
\exp(-5\lambda^{1/3}t/2) \le2\cdot\frac{2}{5}
\exp(-5\lambda^{1/3}t/2),
\end{eqnarray*}
since $\omega-2$ and $\omega^2-2$ each have real part $-5/2$.
Putting all together
%
%
\begin{equation}\label{intbd}
\int_0^t
\exp(-2\lambda^{1/3}s) V'(s) \,ds =\frac{1}{7} -
\frac{1}{3} \exp(-\lambda^{1/3}t)+ \theta_t,
\end{equation}
since $EM_t^2=\exp(-2\lambda^{1/3}t)EY_t^2$, the desired
result follows.
\end{pf}

We use the previous calculation to get bounds for $EA_t^2, EL_t^2$ and
$EX_t^2$, which will be useful later.
\begin{lemma}\label{sqbound}
Let $a(\cdot), l(\cdot)$ and $x(\cdot)$ be as in (\ref{a}). Then
\[
EA_t^2 \le\tfrac{27}{2} a^2(t),\qquad EL_t^2 \le\tfrac{27}{2}
l^2(t),\qquad
EX_t^2 \le\tfrac{27}{2} x^2(t).
\]
\end{lemma}
\begin{pf}
By (\ref{intbd}) we have
%
%
\begin{equation}\label{intbd1}
\int_0^t\exp(-2\lambda^{1/3}s) V'(s) \,ds \le\frac{1}{7}
+ \frac{4}{15} = \frac{43}{105} \le\frac{1}{2}.
\end{equation}
Now using Lemma \ref{renewaleq}
\begin{eqnarray*}
EA_t^2&=&(EA_t)^2+\intt EA_{t-s}^2 \,dF(s),\qquad
EL_t^2=(EL_t)^2+\intt EL_{t-s}^2
\,dF(s),\\
EX_t^2&=&(EX_t)^2+\intt EX_{t-s}^2 \,dF(s).
\end{eqnarray*}
Solving the renewal
equations $EA_t^2=\phi_a*V(t), EL_t^2=\phi_l*V(t)$ and $EX_t^2=\phi_x*V(t)$,
where $V(\cdot)$ is as in Lemma \ref{V}
and $\phi_a(t)=(EA_t)^2, \phi_l(t)=(EL_t)^2$ and $\phi
_x(t)=(EX_t)^2$. A crude upper bound for $\phi_a(t)$ is
$9a^2(t)$. Since $a(t-s)=a(t)\exp(-\lambda^{1/3}s)$,
%
%
\begin{equation}
\label{a2bd}
a^2*V(t)=a^2(t)\biggl[1+\intt\exp(-\lambda^{1/3} s)
V'(s) \,ds\biggr] \le\frac{3a^2(t)}{2}
\end{equation}
by (\ref{intbd1}). Hence $EA_t^2\le9a^2*V(t)\le(27/2)a^2(t)$.

Similarly using the bounds $9l^2(t)$ and $9x^2(t)$ for
$\phi_l(t)$ and $\phi_x(t)$, respectively, and noting that
$l(t-s)/l(t)=x(t-s)/x(t)=\exp(-\lambda^{1/3} s)$, we get
the desired bounds for $EL_t^2$ and $EX_t^2$.
\end{pf}
\begin{lemma} \label{JKbds}
Let $\tilde J_t, \tilde K_t = e^{-\eta t}(X_t + \eta L_t + \eta^2
A_t)$ with $\eta= \omega\lambda^{1/3}$, $\omega^2\lambda^{1/3}$,
respectively. Then $\tilde J_t$ and $\tilde K_t$ are complex
martingales with respect to the filtration $\mathcal{F}_t$, and
\[
E|\tilde J_t|^2, E|\tilde K_t|^2 = \tfrac{1}{6}
\exp(2\lambda^{1/3}t)+\tfrac{1}{2} +\theta_t\qquad \mbox{where }
|\theta_t|\le
\tfrac23 \exp(\lambda^{1/3}t/2),
\]
and hence $E|\tilde J_t|^2, E|\tilde K_t|^2 \le(4/3) \exp
(2\lambda^{1/3}t)$.
\end{lemma}
\begin{pf} Let $h(t,x,\ell,a)=e^{-\eta t}(x+\eta
\ell+\eta^2a)$, and let $\mathcal L$ be the generator of the Markov
process $(t,X_t,L_t,A_t)$. Equation (\ref{infgen}) implies $\mathcal{L}h=0$
when $\eta=\lambda^{1/3}\omega, \lambda^{1/3}\omega^2$, so that
$\tilde J_t$ and $\tilde K_t$ are complex martingales by Dynkin's
formula.\vspace*{1pt}

First we compute $E|J_t|^2$, where
$J_t=\exp(\lambda^{1/3}\omega t) \tilde J_t$. For that we
use Lemma \ref{renewaleq} with $Y_t=J_t$ and $Z_t=\bar J_t$, the
complex conjugate. Since $\tilde J_t$ is a~complex martingale with
$\tilde J_0=1$ and $\omega= (-1+i\sqrt{3})/2$, $E\tilde
J_t=1$ and hence
\[
|EJ_t|^2 = \exp(-\lambda^{1/3} t).
\]
Using Lemma \ref{renewaleq} $E|J_t|^2=|EJ_t|^2 + \intt
E|J_{t-s}|^2 \,dF(s)$. Solving the renewal equation as we have done
twice before
\[
E|J_t|^2 = \exp(-\lambda^{1/3} t) + \int_0^t
\exp\bigl(-\lambda^{1/3}(t-s)\bigr) V'(s) \,ds.
\]
Repeating the first part of the proof for
$K_t=\exp(\lambda^{1/3}\omega^2 t) \tilde K_t$, we see
that $E|K_t|^2$ is also equal to the right-hand side above.

The integral is $\exp(-\lambda^{1/3}t)$ times
\begin{eqnarray*}
&&\frac{1}{3} \int_0^t \exp(\lambda^{1/3}s) \cdot
\lambda^{1/3}[\exp(\lambda^{1/3}s)+\omega
\exp(\lambda^{1/3}\omega s)
+\omega^2 \exp(\lambda^{1/3}\omega^2 s)] \,ds \\
&&\qquad=\frac{1}{3}\biggl[ \frac{1}{1+1}\{
\exp(2\lambda^{1/3}t)-1\}
+\frac{\omega}{\omega+1}\bigl\{\exp\bigl((\omega+1)\lambda
^{1/3} t\bigr)-1\bigr\}\\
&&\qquad\quad\hspace*{108.7pt}{}
+\frac{\omega^2}{\omega^2+1}\bigl\{\exp\bigl((\omega
^2+1)\lambda^{1/3}
t\bigr)-1\bigr\}\biggr].
\end{eqnarray*}
Now using $1= -\omega-\omega^2$ and $\omega^3=1$,
\[
-\frac12-\frac{\omega}{\omega+1}-\frac{\omega^2}{\omega^2+1}
=-\frac
12-\frac{\omega^3+\omega+\omega^3+\omega^2}{\omega^3+\omega
^2+\omega+1}
=-\frac32.
\]
Since $\omega= (-1+i\sqrt{3})/2$ and $\omega^2 =
(-1-i\sqrt{3})/2$, if we take
\[
\theta_t = \frac13\biggl[\frac{\omega}{\omega+1}\exp
\bigl((\omega+1)\lambda^{1/3} t\bigr)
+ \frac{\omega^2}{\omega^2+1}\exp\bigl((\omega^2+1)\lambda^{1/3}
t\bigr) \biggr],
\]
then
\[
3|\theta_t| \le\biggl( \frac{1}{|\omega+1|} +
\frac{1}{|\omega^2+1|} \biggr) \exp(\lambda^{1/3}t/2)
\le
2 \exp(\lambda^{1/3}t/2),
\]
since each of $\omega+1$ and $\omega^2+1$ has real part $1/2$.
Putting all together
%
%
\begin{equation}\label{Jbd}
E|J_t|^2\le\tfrac16 \exp(\lambda^{1/3}t) + \tfrac12
\exp(-\lambda^{1/3}t)
+ \tfrac23 \exp(-\lambda^{1/3}t/2),
\end{equation}
which completes the proof, since
$E|\tilde J_t|^2/E|J_t|^2=\exp(\lambda^{1/3} t)=E|\tilde K_t|^2/E|K_t|^2$.
\end{pf}
\begin{lemma} If $M = \lim_{t\to\infty} M_t$, we have $P(M>0)=1$ and
\[
\exp(-\lambda^{1/3} t) X_t\mbox{, }
\lambda^{1/3}\exp(-\lambda^{1/3} t) L_t\mbox{, }
\lambda^{2/3}\exp(-\lambda^{1/3} t) A_t \to\frac{M}{3}
\]
a.s. and in $L^2$.
\end{lemma}
\begin{pf} $M = \lim_{t\to\infty} M_t$ exists a.s. and in $L^2$,
since $M_t$ is an $L^2$ bounded martingale.
Recall that
\begin{eqnarray*}
I_t &=& X_t + \lambda^{1/3} L_t + \lambda^{2/3} A_t,\\
J_t &=& X_t + \omega\lambda^{1/3} L_t + \omega^2 \lambda^{2/3}
A_t,\\
K_t &=& X_t + \omega^2 \lambda^{1/3} L_t + \omega\lambda^{2/3} A_t.
\end{eqnarray*}
Since $1+\omega+\omega^2=0$ and $\omega^3=1$,
%
%
\begin{eqnarray}\label{lincomb}
3X_t &=& I_t + J_t + K_t,
\nonumber\\
3\lambda^{1/3} L_t &=& I_t + \omega^2 J_t + \omega K_t,
\\
3\lambda^{2/3} A_t &=& I_t + \omega J_t + \omega^2 K_t. \nonumber
\end{eqnarray}
Since $M_t= \exp(-\lambda^{1/3} t) I_t \to M$, it suffices to show
that $\exp(-\lambda^{1/3} t) J_t$ and $\exp(-\lambda^{1/3} t) K_t$
go to 0 a.s. and in $L^2$. We will only prove this for $J_t$, since
the argument for $K_t$ is almost identical. $\tilde J_t$ is a
complex martingale, so $|\tilde J_t|$ is a real submartingale. Using
the $L^2$ maximal inequality, (4.3) in Chapter~4 of \citet{Dur10}
and Lemma \ref{JKbds},
%
%
\begin{equation} \label{L2max}
E\Bigl( \max_{0\le s \le t} |\tilde
J_s|^2 \Bigr) \le4 E|\tilde J_t|^2 \le
\frac{16}{3}\exp(2\lambda^{1/3}t).
\end{equation}
The real part
of $\omega$ is $-1/2$. So writing $\tilde
J_s=\exp(\lambda^{1/3}(1-\omega)s) \cdot\exp(-\lambda^{1/3}s)J_s$,
we see that
%
%
\begin{equation} \label{hammer}
E\Bigl( \max_{u\le s \le t} |\tilde J_s|^2 \Bigr)
\ge\exp(3\lambda^{1/3}u) E\Bigl( {\max_{u \le s \le t}}
|{\exp}(-\lambda^{1/3}s)J_s|^2 \Bigr).
\end{equation}
Combining these bounds with Chebyshev inequality, and taking
$t_n=\break2\lambda^{-1/3}\log n$ for $n=1, 2, \ldots$
%
%
\begin{eqnarray}\label{supJbd}\qquad
P\Bigl( {\max_{t_n \le s \le t_{n+1}}}
|{\exp}(-\lambda^{1/3}s)J_s|^2 \ge\ep\Bigr)
&\le&\ep^{-2} E \Bigl( {\max_{t_n \le s \le t_{n+1}}}
|{\exp}(-\lambda^{1/3}s)J_s|^2 \Bigr) \nonumber
\\
&\le&\frac{16}{3} \ep^{-2} \exp\bigl(\lambda
^{1/3}(2t_{n+1}-3t_n)\bigr)\\
&=& \frac{16}{3}\ep^{-2} \frac{(n+1)^4}{n^6}\nonumber
\end{eqnarray}
for any $\ep>0$. Summing over $n$, and using the
Borel--Cantelli lemma
\[
|{\exp}(-\lambda^{1/3}s)J_s| \to0 \qquad\mbox{a.s.}
\]
To get convergence in $L^2$ we use (\ref{Jbd}).
\[
E|{\exp}(-\lambda^{1/3}t)J_t|^2 \le
\tfrac{4}{3}\exp(-\lambda^{1/3} t) \to0 \qquad\mbox{as }
t\to\infty.
\]

To prove that $P(M>0)=1$ we begin by noting that convergence in $L^2$
implies that $P(M>0)>0$. Every time a new balloon is born it has
positive probability of starting a process with a positive limit,
so this will happen eventually and $P(M>0)=1$.
\end{pf}

\section{Hitting times for $\mathcal{A}_t$ and $\mathcal{C}_t$}\label{sec3}

Recall that $\sigma(\ep) = \inf\{ t \dvtx A_t \ge\ep N^2 \}$ and
$\tau(\ep)=\inf\{t\dvtx C_t\ge\ep N^2\}$. Also recall the definitions
of $a(\cdot), l(\cdot), x(\cdot)$ and $S(\cdot)$ from (\ref{a}) and
(\ref{S}). Note that $a(S(\ep)) = \ep N^2$ and $A_t/a(t), L_t/l(t),
X_t/\allowbreak x(t) \to M$ a.s. by Theorem \ref{th1}. We begin by estimating
the difference between~$M$ and each of $A_t/a(t), L_t/l(t)$ and
$X_t/x(t)$.
\begin{lemma}\label{supbound}
For any $\gamma, u>0$
\[
P\Bigl({\sup_{t\ge u}} |A_t/a(t)-M| \ge\gamma^2\Bigr)
\le C\gamma^{-4}\exp(-\lambda^{1/3} u)
\]
for some constant $C$. The same bound holds for $ P({\sup_{t\ge u}}
|L_t/l(t)-M| \ge\gamma^2)$ and $P({\sup_{t\ge u}}
|X_t/x(t)-M| \ge\gamma^2)$.
\end{lemma}
\begin{pf}
Using (\ref{lincomb}) $A_t/a(t)=M_t+\omega
\exp(-\lambda^{1/3}t) J_t+\omega^2
\exp(-\lambda^{1/3}t) K_t$. For $0<u\le t$ the triangle
inequality
implies
%
%
\begin{equation}\label{bd1}
\quad|A_t/a(t) - M| \le|M_t-M| +
|{\exp}(-\lambda^{1/3}t) J_t| +
|{\exp}(-\lambda^{1/3}t) K_t|.
\end{equation}
Taking the supremum over $t$,
%
%
\begin{eqnarray}\label{supbd}
&&P\Bigl({\sup_{t\ge u} }|A_t/a(t)-M| \ge\gamma^2\Bigr)\nonumber\\
&&\qquad
\le P\Bigl({\sup_{t\ge u} }|M_t-M| \ge\gamma^2/3\Bigr)
+ P\Bigl({\sup_{t\ge u} }|{\exp}(-\lambda^{1/3}t)
J_t| \ge\gamma^2/3\Bigr)\\
&&\qquad\quad{} + P\Bigl({\sup_{t\ge u}} |{\exp}(-\lambda^{1/3}t)
K_t| \ge\gamma^2/3\Bigr).\nonumber
\end{eqnarray}
To bound the first term in the right-hand side of (\ref{supbd}) we
note that
\[
E\Bigl({\sup_{t\ge u}}|M_t-M|^2\Bigr) = \lim_{U\to\infty}
E\Bigl({\max_{u\le t\le U}} |M_t-M|^2\Bigr).
\]
Using triangle inequality $|M_t-M|\le|M_t-M_u|+|M_u-M|$. Taking
supremum over $t\in[u,U]$
and using the inequality $(a+b)^2\le2(a^2+b^2)$,
\[
E\Bigl({\max_{u\le t\le U}} |M_t-M|^2\Bigr) \le
2\Bigl(E\Bigl({\max_{u\le t\le U}} |M_t-M_u|^2
\Bigr)+E|M_u-M|^2\Bigr).
\]
Using the $L^2$
maximal inequality, (4.3) in Chapter 4 of \citet{Dur10} and
orthogonality of martingale increments,
\[
E\Bigl({\max_{u\le t\le U} }|M_t-M_u|^2\Bigr) \le
4E(M_U-M_u)^2=4(EM_U^2-EM_u^2).
\]
Since the martingale $M_t$ converges to $M$ in $L^2$,
$EM^2=\lim_{t\to\infty} EM_t^2=8/7$. Then using orthogonality of
martingale increments and Lemma \ref{mart},
\[
E(M_u-M)^2 = EM^2 - EM_u^2 \le\exp(-\lambda^{1/3}u).
\]
Combining the last four bounds with Lemma \ref{mart}, and using
Chebyshev inequality
%
%
\begin{equation}\label{supbd1}
P\Bigl({\sup_{t\ge u}} |M_t-M| \ge\gamma^2/3\Bigr) \le
9\gamma^{-4}\cdot10\exp(-\lambda^{1/3} u).
\end{equation}

To bound the second term in the right-hand side of (\ref{supbd}) we
take $t_n=u+2\lambda^{-1/3}\log n$ for $n=1, 2, \ldots$ and use an
argument similar to the one leading to (\ref{supJbd}) together with
Chebyshev inequality to get
%
%
\begin{eqnarray} \label{supbd2}
P\Bigl({\sup_{t\ge u}}
|{\exp}(-\lambda^{1/3}t) J_t| \ge \gamma
^2/3\Bigr)
&\le&\sum_{n=1}^\infty P\Bigl({\max_{t_n\le t\le t_{n+1}}}
|{\exp}(-\lambda^{1/3}t) J_t| \ge\gamma
^2/3\Bigr)
\nonumber\\
&\le& 9\gamma^{-4} \sum_{n=1}^\infty E\Bigl({\max_{t_n\le t\le t_{n+1}}}
|{\exp}(-\lambda^{1/3}t) J_t|\Bigr)^2
\nonumber\\[-8pt]\\[-8pt]
&\le& 9 \cdot\frac{16}{3}\gamma^{-4} \sum_{n=1}^\infty
\exp\bigl(\lambda^{1/3}(2t_{n+1}-3t_n)\bigr)
\nonumber\\
&=& 48 \gamma^{-4} \exp(-\lambda^{1/3}u)
\sum_{n=1}^\infty\frac{(n+1)^4}{n^6}.\nonumber
\end{eqnarray}
Repeating the previous argument for the third term in the right-hand
side of (\ref{supbd})
we get the same upper bound as in (\ref{supbd2}). Combining
(\ref{supbd}), (\ref{supbd1}) and~(\ref{supbd2})
we get the desired bound for $A_t/a(t)$.

The bound in (\ref{bd1}) also works for both $L_t/l(t)$ and
$X_t/x(t)$, since
using~(\ref{lincomb})
\begin{eqnarray*}
L_t/l(t) &=& M_t+\omega^2\exp(-\lambda^{1/3}t)J_t + \omega\exp
(-\lambda^{1/3}t)K_t,\\
X_t/x(t) &=& M_t+ \exp(-\lambda^{1/3}t)J_t+
\exp(-\lambda^{1/3}t)K_t,
\end{eqnarray*}
and so the assertion of this lemma holds if $A_t/a(t)$ is replabed
by $L_t/l(t)$ or~$X_t/x(t)$.
\end{pf}

We now use Lemma \ref{supbound} to study the limiting behavior of
$\sigma(\ep)$.
\begin{lemma}\label{ALXbd}
Let $W_\ep=S(\ep/M)$, where $S(\cdot)$ is as in (\ref{S}) and $M$
is the limit random variable in Theorem \ref{th1}.
Then for any $\eta>0$
\begin{eqnarray*}
\lim_{N\to\infty} P(|A_{W_\ep}-\ep N^2|>\eta N^2)
&=& \lim_{N\to\infty} P(|L_{W_\ep}-\ep N^{2-\alpha/3}|>\eta
N^{2-\alpha/3})\\
&=& \lim_{N\to\infty} P(|X_{W_\ep}-\ep N^{2-2\alpha/3}|>\eta
N^{2-2\alpha/3})\\
&=&0.
\end{eqnarray*}
\end{lemma}
\begin{pf}
Since $P(M>0)=1$, given $\theta>0$, we can choose $\gamma=\gamma
(\theta)>0$ so
that $\gamma<\eta/\ep$ and
%
%
\begin{equation}\label{Mnot0}
P(M<\gamma)<\theta.
\end{equation}
Using Lemma \ref{supbound} we can choose a constant
$b=b(\gamma,\theta)$ such that
\[
P\Bigl( {\sup_{t\ge b N^{\alpha/3}}} |A_t/a(t) - M | > \gamma^2
\Bigr) <
\theta.
\]
Combining with (\ref{Mnot0})
\[
P\Bigl( {\sup_{t\ge bN^{\alpha/3}}} |A_t/a(t) - M | > \gamma M
\Bigr) <
2\theta.
\]
Since $a(W_\ep)=\ep N^2/M$, by the choices of $\gamma$ and $b$,
\begin{eqnarray*}
P( |A_{W_\ep} - \ep N^2| \ge\eta N^2) &\le& P( |A_{W_\ep} - \ep
N^2| \ge\ep\gamma N^2) \\
&=& P\bigl(|A_{W_\ep}/a(W_\ep)-M|\ge\gamma M\bigr) \\
&<& 2\theta
+P(W_\ep<bN^{\alpha/3}).
\end{eqnarray*}
By the definition of $S(\cdot)$,
\[
P( W_\ep< b N^{\alpha/3}) = P\biggl(M > \frac{3\ep}{b}
N^{2-2\alpha/3}\biggr) \to0
\]
as $N\to\infty$, and so $\limsup_{N\to\infty} P(|A_{W_\ep
}-\ep N^2|>\eta N^2) \le2\theta$.
Since $\theta>0$ is arbitrary, we have shown that
\[
\lim_{N\to\infty} P( |A_{W_\ep}-\ep N^2| \ge\eta N^2
) =
0.
\]
Repeating\vspace*{1pt} the argument for $L_{W_\ep}$ and $X_{W_\ep}$, and noting
that $l(W_\ep)=\break\ep N^{2-\alpha/3}/M$ and $x(W_\ep)=\ep
N^{2-2\alpha/3}/M$, we get the other two assertions.~%
\end{pf}

As a corollary of Lemma \ref{ALXbd} we get the first conclusion of
Theorem \ref{th2}.
\begin{corollary} \label{th2part1}
As $N\to\infty$, $N^{-\alpha/3} (\sigma(\ep)-S(\ep)) \eqp-\log(M)$.
\end{corollary}
\begin{pf}
For any $\eta>0$ choose $\gamma>0$ so that $\log(1+\gamma)<\eta$
and $\log(1-\gamma)>-\eta$. Let $W_\ep$ be as in Lemma \ref{ALXbd}.
Clearly $W_{(1+\gamma)\ep}=S(\ep)+N^{\alpha/3}[\log(1+\gamma
)-\log M]$ and $W_{(1-\gamma)\ep}=S(\ep)+N^{\alpha/3}[\log
(1-\gamma)-\log M]$.
Using Lemma~\ref{ALXbd}
\begin{eqnarray*}
&& P\bigl[N^{-\alpha/3} \bigl(\sigma(\ep)-S(\ep)\bigr)>-\log M+\eta\bigr]
\\
&&\qquad \le P\bigl(\sigma(\ep) >W_{(1+\gamma)\ep}\bigr)
=P(A_{W_{(1+\gamma)\ep}}<\ep N^2)\to0,\\
&& P\bigl[N^{-\alpha/3} \bigl(\sigma(\ep)-S(\ep)\bigr)<-\log M-\eta\bigr]
\\
&&\qquad \le P\bigl(\sigma(\ep) < W_{(1-\gamma)\ep}\bigr)
=P(A_{W_{(1-\gamma)\ep}} > \ep N^2)\to0
\end{eqnarray*}
as $N\to\infty$, and the proof is complete.
\end{pf}

The second conclusion in Theorem \ref{th2} follows from $C_t \le A_t$. To
get the third we have to wait till Lemma \ref{tausigma}. First we
need to show that when~$A_t/N^2$ is small, $C_t/N^2$ is not very
much smaller. To prepare for that we need
the following result.
\begin{lemma}\label{renewalineq}
Let $F(t)=\lambda t^3/3!$. If $u(\cdot)$ and $\beta(\cdot)$ are
functions such that $u(t) \le\beta(t)+\intt u(t-s) \,dF(s)$ for all
$t\ge0$, then
\[
u(t) \le\beta* V(t) = \beta(t)+\intt\beta(t-s) \,dV(s),
\]
where $V(\cdot)$ is as in Lemma \ref{V}.
\end{lemma}
\begin{pf}
Define $\tilde\beta(t)\equiv\beta(t)+\intt u(t-s) \,dF(s)-u(t)$.
So $\tilde\beta(t) \ge0$ for all $t\ge0$. If $\hat\beta(t)\equiv
\beta(t)-\tilde\beta(t)$, then
\[
u(t)=\hat\beta(t)+\intt u(t-s) \,dF(s).
\]
Solving the renewal equation we get
$u(t)=\hat\beta* V(t)$, where
$V(\cdot)$ is as in Lem\-ma~\ref{V}. Since $\hat\beta(t)\le\beta
(t)$ for
all $t\ge0$, we get the result.
\end{pf}

We now apply Lemma \ref{renewalineq} to estimate the difference
between $EA_t$ and~$EC_t$.

\begin{lemma}\label{compare1}
For any $t\ge0$ and $a(\cdot)$ as in (\ref{a}),
\[
EC_t\ge EA_t - \frac{11 a^2(t)}{ N^2}.
\]
\end{lemma}
\begin{pf}
In either of our processes, if a center is born at time $s$, then
the radius of the corresponding disk at time $t>s$ will be
$(t-s)/\sqrt{2\pi}$. Thus $x$ will be covered at time $t$ if and
only if there is a center in\vadjust{\goodbreak} the space--time cone
%
%
\begin{equation}\label{cone}
K_{x,t}\equiv\bigl\{(y,s)\in\Gamma(N) \times[0,t]\dvtx|y-x| \le
(t-s)/\sqrt{2\pi}\bigr\}.
\end{equation}
If $0=s_0, s_1, s_2,\ldots$ are the
birth times of new centers in $\mathcal{C}_t$, then
\[
P( x \notin\mathcal{C}_t | s_0, s_1, s_2, \ldots) = \prod_{i\dvtx
s_i\le
t} \biggl[1-\frac{(t-s_i)^2}{2N^2}\biggr] \le\exp\biggl[-\sum_{i\dvtx
s_i\le t} \frac{(t-s_i)^2}{2N^2}\biggr],
\]
since $1-x\le e^{-x}$. Let $q(t)\equiv P(x \notin\mathcal{C}_t)$,
which does not depend on $x$, since we have a random chosen starting
point. Recall that $\tilde X_t$ is the number of centers born by
time $t$ in $\mathcal{C}_t$. Using the last inequality
\[
q(t) \le E\exp\biggl[-\int_0^t \frac{(t-s)^2}{2N^2} \,d\tilde
X_s\biggr]
\]
and $E C_t= N^2(1- q(t))$. Integrating $e^{-y} \ge1-y$ gives
$1-e^{-x}\ge x-x^2/2$ for $x\ge0$. So
%
%
\begin{eqnarray}
\label{eq7}
E C_t & \ge & N^2E\biggl[1-\exp\biggl(-\int_0^t \frac
{(t-s)^2}{2N^2} \,d\tilde X_s\biggr)\biggr]\nonumber\\[-8pt]\\[-8pt]
& \ge & N^2 E\biggl[\int_0^t \frac{(t-s)^2}{2N^2} \,d\tilde X_s
-\frac12\biggl(\int_0^t \frac{(t-s)^2}{2N^2} \,d\tilde X_s
\biggr)^2\biggr]. \nonumber
\end{eqnarray}
For the first term on the right we use $E\tilde X_t=1+\lambda\intt
EC_s \,ds$. For the second term on the right, we use the coupling
between $\mathcal{C}_t$ and $\mathcal{A}_t$ described in the \hyperref
[intro]{Introduction},
see (\ref{couple}), so that we have $\int_0^t (t-s)^2 \,d\tilde X_s
\le\int_0^t (t-s)^2 \,dX_s$. Combining these two facts
%
%
\begin{eqnarray}\label{eq6}
EC_t & \ge & \frac{t^2}{2} + \int_0^t \frac{(t-s)^2}{2} \lambda EC_s \,ds
-\frac{1}{2N^2} E\biggl[\int_0^t \frac{(t-s)^2}{2}\,dX_s\biggr]^2
\nonumber\\[-8pt]\\[-8pt]
&=& \frac{t^2}{2} + \intt\frac{(t-s)^2}{2} \lambda EC_s \,ds-\frac
{EA_t^2}{2N^2}.\nonumber
\end{eqnarray}
The last equality follows from (\ref{LA}), as does the next equation
for $EA_t$:
%
%
\begin{equation}\label{eq13}
EA_t=\frac{t^2}{2} +\int_0^t
\frac{(t-s)^2}{2} V'(s) \,ds = \frac{t^2}{2} + \int_0^t
\frac{(t-s)^2}{2} \lambda EA_s \,ds.
\end{equation}
Here $V(\cdot)$ is as in
Lemma \ref{V} and $EA_t=V'(t)/\lambda$ by Lemma \ref{XLAlem}.
Combining~(\ref{eq6}) and (\ref{eq13}), if $u(t)\equiv EA_t-EC_t$,
and $F(s)=\lambda s^3/3!$, then
\[
u(t) \le\frac{EA_t^2}{2N^2} + \int_0^t \frac{(t-s)^2}{2} \lambda
u(s) \,ds = \frac{EA_t^2}{2N^2} + \int_0^t u(t-r) \,dF(r),
\]
where the last step is obtained by changing variables $s \mapsto
t-r$. If $\beta(t) = EA_t^2/2N^2$, then by Lemma \ref{sqbound}
$\beta(t)\le27a^2(t)/4N^2$, and using\vadjust{\goodbreak} Lemma \ref{renewalineq} and
(\ref{a2bd})
\[
u(t) \le\beta* V(t)\le\frac{27}{4N^2} (a^2)*V(t) \le
\frac{27}{4N^2} \frac32 a^2(t),
\]
which gives the result, since $81/8 \le11$.
\end{pf}

To complete the proof of Theorem \ref{th2} it remains to show the
third conclusion of it, which we separate as the following lemma and
prove it using
Lemma \ref{compare1}.
\begin{lemma} \label{tausigma}
For any $\gamma>0$
\[
\limsup_{N\to\infty} P\bigl( \tau(\ep) > \sigma\bigl((1+\gamma)\ep\bigr) \bigr) \le
P\bigl( M \le(1+\gamma)\ep^{1/3} \bigr) + 11
\frac{\ep^{1/3}}{\gamma}.
\]
\end{lemma}
\begin{pf}
Let $U=\sigma((1+\gamma)\ep)$ and $T=S(\ep^{2/3})$, where $S(\cdot)$
is as in (\ref{S}). Now
\[
S(\ep^{2/3}) - S\bigl((1+\gamma)\ep\bigr) = N^{\alpha/3} \bigl[ - \tfrac{1}{3}
\log(\ep) - \log(1+\gamma) \bigr].
\]
It follows from Corollary \ref{th2part1} that
\begin{eqnarray*}
\limsup_{N\to\infty}
P( U \ge T ) &\le& P\biggl( - \log(M) \ge- \frac{1}{3} \log(\ep) - \log(1+\gamma)
\biggr) \\
&=& P\bigl( M \le(1+\gamma) \ep^{1/3} \bigr).
\end{eqnarray*}
Using Markov's inequality, Lemma \ref{compare1}, and $a(T) =
\ep^{2/3}N^2$,
%
%
\begin{equation}\label{b3}\quad
P(|A_{T}-C_{T}|>\gamma\ep N^2) \le
\frac{E(A_{T}-C_{T})}{\gamma\ep N^2} \le\frac{6 (a(T))^2}{\gamma
\ep
N^4} \le11 \cdot\frac{\ep^{1/3}}{\gamma}.
\end{equation}
Using these two
bounds and the fact that $|A_t-C_t|$ is nondecreasing in $t$, we get
\begin{eqnarray*}
&&\limsup_{N\to\infty}
P\bigl[\tau(\ep)>\sigma\bigl((1+\gamma)\ep\bigr)\bigr]\\
&&\qquad= \limsup_{N\to\infty} P[|A_U-C_U| >\gamma\ep
N^2]\\
&&\qquad \le\limsup_{N\to\infty} P( U \ge T ) + \limsup_{N\to\infty}
P[|A_U-C_U| >\gamma\ep N^2, U < T]\\
&&\qquad \le\limsup_{N\to\infty} P( U \ge T ) + P(
|A_{T}-C_{T}|>\gamma\ep N^2),
\end{eqnarray*}
which completes the proof.
\end{pf}

\section{Limiting behavior of $\mathcal{C}_t$}\label{sec4}

Let $\mathcal{C}_{s,t}^0$ be the set of points covered in $\mathcal{C}_t$
at time $t$ by the balloons born before time $s$.
If we number the generations of centers in $\mathcal{C}_t$
starting with those existing at time $s$ as $\mathcal{C}_t$-centers of
generation~0, then $\mathcal{C}_{s,t}^0$ is the set of points covered
at time $t$
by the generation~0 centers of\vadjust{\goodbreak} $\mathcal{C}_t$. Let $\mathcal
{C}^1_{s,t}$ be the set of
points, which are either
in $\mathcal{C}^0_{s,t}$, or are covered at time $t$ by a balloon born
from this area. This
is the set of points covered by $\mathcal{C}_t$-centers of generations
$\le1$ at time~$t$, ignoring births from
$\mathcal{C}^1_{s,t} \setminus\mathcal{C}^0_{s,t}$, which are
second generation centers. Continuing by induction, we let $\mathcal{
C}^k_{s,t}$ be the set of points and $C_{s,t}^k=|\mathcal
{C}_{s,t}^k|$ be the total area covered by $\mathcal
{C}_t$-centers of generations $0\le j \le k$ at time $t$. Similarly
$A_{s,t}^k$ denotes the total area of the balloons in $\mathcal{A}_t$
of generations $j\in\{0, 1, \ldots, k\}$ at time~$t$, where
generation 0 centers are those existing at time $s$.

Recall the following definitions from (\ref{a}), (\ref{S}),
(\ref{R}) and (\ref{psiWI}).
\begin{eqnarray*}
a(t) &=& (1/3)N^{2\alpha/3}\exp(N^{-\alpha/3} t),\\
S(\ep) &=& N^{\alpha/3}[(2-2\alpha/3)\log N + \log(3\ep)], \\
R &=& N^{\alpha/3}[(2-2\alpha/3)\log N - \log(M)],
\end{eqnarray*}
where $M$ is the limit random variable in Theorem \ref{th1}, and for
$\log(3\ep) \le t$,
\[
\psi(t)\equiv R+N^{\alpha/3}t,\qquad W\equiv\psi(\log(3\ep))\quad
\mbox{and}\quad I_{\ep,t}=[\log(3\ep) , t].
\]
Note that $\psi(t) \le0$ only if $M\ge N^{2-2\alpha/3}t$.\vspace*{1pt}

Obviously $C_{s,t}^0 \le A_{s,t}^0$.
For the other direction we have the following lemma.
\begin{lemma}\label{compare2}
For any $0<s<t$,
\[
EC_{s,t}^0 \ge EA_{s,t}^0-\frac{a^2(s)}{N^2}p\bigl((t-s)\lambda
^{1/3}\bigr),
\]
where for some positive constants $c_1, c_2$ and $c_4$,
%
%
\begin{equation} \label{pxdef}
p(x)=c_1+c_2x^2/2!+c_4x^4/4!.
\end{equation}
\end{lemma}
\begin{pf} By the definition of $A_{s,t}^0$,
%
%
\begin{equation}\label{Ast}
A_{s,t}^0=\int_0^s \frac{(t-r)^2}{2}
\,dX_r =
\frac{(t-s)^2}{2}X_s+(t-s)L_s+A_s.
\end{equation}
For the second equality we have written
$(t-r)^2=(t-s)^2+2(t-s)(s-r)+(s-r)^2$ and used (\ref{LA}).
As in Lemma \ref{compare1}, a point $x$ is not covered by time~$t$
by the balloons born before time $s$, if and only if no center is born
in the
truncated space--time cone
\[
K_{x,s,t} \equiv\bigl\{(y,r)\in\Gamma(N)\times[0,s]\dvtx
|y-x| \le(t-r)/\sqrt{2\pi}\bigr\}.
\]
So using arguments similar to the ones for (\ref{eq7}) and $1-e^{-x}
\ge x - x^2/2$,
\begin{eqnarray*}
EC_{s,t}^0 &\ge& N^2E\biggl[1-\exp\biggl(-\int_0^s \frac{(t-r)^2}{2N^2}
\,d\tilde X_r\biggr)\biggr]\\
&\ge& N^2\biggl[E\int_0^s \frac{(t-r)^2}{2N^2}
\,d\tilde X_r-\frac12 E\biggl(\int_0^s \frac{(t-r)^2}{2N^2}
\,d\tilde X_r\biggr)^2\biggr].
\end{eqnarray*}
For the first term on the right, we use $E\tilde X_t=1+\lambda\intt
EC_s \,ds$. For the second term on the right, we
use the coupling between $\mathcal{C}_t$ and $\mathcal{A}_t$
described in
the \hyperref[intro]{Introduction}, see (\ref{couple}), to conclude that
\[
\int_0^s (t-r)^2 \,d\tilde X_r \le\int_0^s (t-r)^2 \,dX_r=2A_{s,t}^0.
\]
Combining these two facts, using the first equality in (\ref{Ast}),
$EX_t=1+\lambda\intt EA_s \,ds$, and Lemma \ref{compare1},
%
%
\begin{eqnarray}\label{eq1}\quad
EC_{s,t}^0 & \ge & \frac{t^2}{2} + \int_0^s \frac{(t-r)^2}{2} \lambda
EC_r \,dr
-\frac{E(A_{s,t}^0)^2}{2N^2}
\nonumber\\[2pt]
& \ge & \frac{t^2}{2} + \int_0^s \frac{(t-r)^2}{2} \lambda EA_r \,dr
- 11\int_0^s \frac{(t-r)^2}{2}\frac{\lambda a^2(r)}{N^2} \,dr-\frac
{E(A_{s,t}^0)^2}{2N^2}
\\[2pt]
&=& EA_{s,t}^0 - 11\int_0^s \frac{(t-r)^2}{2}\frac{\lambda a^2(r)}{N^2}
\,dr-\frac{E(A_{s,t}^0)^2}{2N^2}.\nonumber
\end{eqnarray}
To estimate the second term in the right-hand side of (\ref{eq1}), we
write
\[
(t-r)^2/2=(t-s)^2/2+(t-s)(s-r)+(s-r)^2/2,
\]
change variables $r = s-q$, and note $a(s-q)=a(s)\exp(-\lambda
^{1/3}q)$, to get
%
%
\begin{eqnarray}\label{2nd}
&&\int_0^s\frac{(t-r)^2}{2} \lambda a^2(r) \,dr\nonumber\\[2pt]
&&\qquad= a^2(s) \biggl[\frac{(t-s)^2}{2}\lambda^{2/3}
\int_0^s\lambda^{1/3}\exp(-2\lambda^{1/3}q) \,dq
\nonumber\\[2pt]
&&\qquad\quad\hspace*{28.1pt}{} +(t-s)\lambda^{1/3} \int_0^s\lambda^{2/3} q
\exp(-2\lambda^{1/3}q) \,dq\\[2pt]
&&\qquad\quad\hspace*{84pt}{}+ \int_0^s\lambda\frac{q^2}{2} \exp(-2\lambda^{1/3}q
) \,dq\biggr]
\nonumber\\[2pt]
&&\qquad\le\frac{a^2(s)}{2}\biggl[\frac{(t-s)^2}{2}\lambda
^{2/3}+(t-s)\lambda^{1/3}+1\biggr].\nonumber
\end{eqnarray}
For the last inequality we have used
\[
\int_0^s r^k\exp(-\mu r) \,dr \le\int_0^\infty r^k\exp(-\mu r)
\,dr = \frac{k!}{\mu^{k+1}}.
\]
To estimate the third term in the right-hand side of (\ref{eq1})
we use (\ref{Ast}) to get
\[
E[(A_{s,t}^0)^2] \le
3[EX_s^2(t-s)^4/4+EL_s^2(t-s)^2+EA_s^2].
\]
Applying Lemma \ref{sqbound} and using the fact that
$a(s)=\lambda^{-1/3}l(s)=\lambda^{-2/3}x(s)$,
%
%
\begin{eqnarray}\label{3rd}
E[(A_{s,t}^0)^2] & \le & 3 \cdot\frac{27}{2}
\biggl[x^2(s)\frac{(t-s)^4}{4} + l^2(s) (t-s)^2 + a^2(s)\biggr]
\nonumber\\[-9pt]\\[-9pt]
& \le & 243 a^2(s)
\biggl[\frac{(t-s)^4}{4!}\lambda^{4/3} + \frac{(t-s)^2}{2!}\lambda
^{2/3}+1\biggr].\nonumber
\end{eqnarray}
Combining (\ref{eq1}), (\ref{2nd}) and (\ref{3rd}) we get the result.
\end{pf}

To show uniform convergence of $C_{W,\psi(\cdot)}^k$ to $C_{\psi
(\cdot)}$, we also need to bound the difference $A_t$ and $A_{s,t}^k$ for
suitable choices of $s$ and $t$.
\begin{lemma}\label{Abound}
If $T=S(\ep^{2/3})$, where $S(\cdot)$ is as in (\ref{S}), then for
any $t>0$
\[
EA_{T+tN^{\alpha/3}}-EA_{T,T+tN^{\alpha/3}}^k \le
3\ep^{2/3}N^2\sum_{j=k+1}^\infty\frac{t^j}{j!}.
\]
\end{lemma}
\begin{pf} By (\ref{Ast})
$EA_{s,t}^0=EA_s+EL_s (t-s)+EX_s (t-s)^2/2$. If $X_{s,t}^k$
and~$L_{s,t}^k$ denote the number of centers and sum of radii of
all the balloons in $\mathcal{A}_t$ of generations $j\in\{1, 2,
\ldots,
k\}$ at time $t$, where generation 0 centers are those which
are born before time $s$, then for $t>s$,
\[
\frac{d}{dt} EX_{s,t}^1 = N^{-\alpha}EA_{s,t}^0,\qquad
\frac{d}{dt} EL_{s,t}^1 = EX_{s,t}^1,\qquad
\frac{d}{dt} EA_{s,t}^1 = EL_{s,t}^1.
\]
Integrating over $[s,t]$ and using (\ref{Ast}) we have
\begin{eqnarray*}
EX_{s,t}^1 &=& N^{-\alpha}\biggl[(t-s)EA_s+\frac{(t-s)^2}{ 2!}
EL_s+\frac{(t-s)^3}{3!}EX_s\biggr],\\[-2pt]
EL_{s,t}^1 &=& N^{-\alpha}\biggl[\frac{(t-s)^2}{2!}EA_s+\frac{(t-s)^3}{3!}
EL_s+\frac{(t-s)^4}{4!}EX_s\biggr],\\[-2pt]
EA_{s,t}^1 &=& N^{-\alpha}\biggl[\frac{(t-s)^3}{3!}EA_s+\frac{(t-s)^4}{4!}
EL_s+\frac{(t-s)^5}{5!}EX_s\biggr].
\end{eqnarray*}
Turning to other generations, for $k\ge2$ and $t>s$,
\begin{eqnarray*}
\frac{d}{dt} (EX_{s,t}^k-EX_{s,t}^{k-1}) &=& N^{-\alpha
}(EA_{s,t}^{k-1}-EA_{s,t}^{k-2}),\\[-2pt]
\frac{d}{dt} (EL_{s,t}^k-EL_{s,t}^{k-1}) &=&
(EX_{s,t}^k-EX_{s,t}^{k-1}),\\[-2pt]
\frac{d}{dt} (EA_{s,t}^k-EA_{s,t}^{k-1}) &=&
(EL_{s,t}^k-EL_{s,t}^{k-1}),
\end{eqnarray*}
and using induction on $k$ we have
\[
EA_{s,t}^k=\sum_{j=0}^k N^{-\alpha j}\biggl[\frac{(t-s)^{3j}}{(3j)!}
EA_s + \frac{(t-s)^{3j+1}}{(3j+1)!} EL_s
+\frac{(t-s)^{3j+2}}{(3j+2)!} EX_s\biggr].\vadjust{\goodbreak}
\]
Since $A_{s,t}^k \uparrow A_t$ for any $s<t$, $EA_t=\lim_{k\to\infty}
EA_{s,t}^k$ by Monotone Convergence theorem. Replacing $s$ by $T$ and
$t$ by $T+t N^{\alpha/3}$,
%
%
\begin{eqnarray}\label{eq5}\qquad
&& EA_{T+t N^{\alpha/3}}-EA_{T,T+t N^{\alpha/3}}^k
\nonumber\\[-8pt]\\[-8pt]
&&\qquad=\sum_{j=k+1}^\infty\biggl[\frac{t^{3j}}{(3j)!} EA_T
+\frac{t^{3j+1}}{(3j+1)!} N^{\alpha/3} EL_T+ \frac
{t^{3j+2}}{(3j+2)!} N^{2\alpha/3} EX_T\biggr].
\nonumber
\end{eqnarray}
Using the fact that $EA_T+ N^{\alpha/3}EL_T+N^{2\alpha
/3}EX_T-3a(T)=0$ and
$a(T)=\ep^{2/3}N^2$, the right-hand side of (\ref{eq5}) is
$\le3\ep^{2/3}N^2\sum_{j=k+1}^\infty t^j/j!$, which completes the proof.
\end{pf}

Recall the definitions of $\psi(\cdot), W$ and $I_{\ep,t}$ from the
displays before Lem\-ma~\ref{compare2} and that for $\log(3\ep)\le t$,
%
%
\begin{equation}\label{gdef2}
g_0(t)=\ep\biggl[1+\bigl(t-\log(3\ep)\bigr)+\frac{(t-\log(3\ep))^2}{2}\biggr].
\end{equation}
\begin{lemma}\label{B0bounds}
For any $t<\infty$, there is an $\ep_0=\ep_0(t)>0$ so that for $0<
\ep< \ep_0$,
\begin{eqnarray*}
\lim_{N\to\infty} P\Bigl( \sup_{s\in I_{\ep,t}}
\bigl|N^{-2} A^0_{W,\psi(s)} -g_0(s) \bigr|>\eta\Bigr)
&=& 0 \qquad\mbox{for any } \eta>0,\\
P\Bigl( \inf_{s\in I_{\ep,t}}
N^{-2} \bigl( C^0_{W,\psi(s)} - A^0_{W,\psi(s)} \bigr) < -\ep
^{7/6} \Bigr)
& \le & P(M<\ep^{1/3})+\ep^{1/12}.
\end{eqnarray*}
\end{lemma}
\begin{pf} To prove the first result we use (\ref{Ast}) to conclude
\[
A^0_{W,\psi(t)}=\frac{(t-\log(3\ep))^2}{2}N^{2\alpha/3}X_W+\bigl(t-\log
(3\ep)\bigr)N^{\alpha/3}L_W+A_W.
\]
Applying Lemma \ref{ALXbd}
\begin{eqnarray*}
&&\lim_{N\to\infty}  P\Bigl(\sup_{s\in I_{\ep,t}}
\bigl|N^{-2} A^0_{W,\psi(s)} -g_0(s)\bigr|>\eta\Bigr) \\
&&\qquad \le\lim_{N\to\infty} P\biggl( \bigl|N^{-(2-2\alpha/3)}X_W-\ep\bigr| >
\frac{2\eta}{3(t-\log(3\ep))^2} \biggr)\\
&&\qquad\quad{} + \lim_{N\to\infty} P\biggl( \bigl|N^{-(2-\alpha/3)}L_W-\ep\bigr| >
\frac{\eta}{3(t-\log(3\ep))} \biggr) \\
&&\qquad\quad{} + \lim_{N\to\infty} P\biggl(|N^{-2}A_W-\ep| >
\frac{\eta}{3}\biggr) =0.
\end{eqnarray*}

Let $\ep_0=\ep_0(t)$ be such that $\ep_0^{1/12}p(t-\log(3\ep)) \le
1$, where $p(\cdot)$ is the polynomial in (\ref{pxdef}).
Let $T=S(\ep^{2/3})$, where $S(\cdot)$ is defined in (\ref{S}), and
$T' = T+(t-\log(3\ep))N^{\alpha/3}$. Using the fact that\vadjust{\goodbreak}
$A^0_{s,s+t}-C^0_{s,s+t}$ is nondecreasing in~$s$, Markov's
inequality, and then Lemma \ref{compare2} we see that
\begin{eqnarray*}
&& P\Bigl(\sup_{s\in I_{\ep,t}} \bigl|A^0_{W,\psi(s)}-C^0_{W,\psi
(s)}\bigr| > \ep^{7/6} N^2, W \le T\Bigr)\\
&&\qquad \le P(|A^0_{T,T'} - C^0_{T,T'}| >\ep^{7/6} N^2)
\le\frac{E|A^0_{T,T'} - C^0_{T,T'}|}{\ep^{7/6} N^2} \\
&&\qquad \le\frac{a^2(T) p(t-\log(3\ep))}{\ep^{7/6} N^4}.
\end{eqnarray*}
Noting that $P(W >T)=P(M<\ep^{1/3}), a(T)=\ep^{2/3}N^2$ and
$\ep^{1/12}p(t-\break\log(3\ep))\hspace*{-0.2pt}<1$ for $\ep<\ep_0$ we have
\[
P\Bigl(\sup_{s\in I_{\ep,t}} \bigl|A_{W,\psi(s)}-C_{W,\psi
(s)}\bigr|>\ep^{7/6}
N^2\Bigr) \le P(M<\ep^{1/3}) + \ep^{1/12},
\]
which completes the proof.
\end{pf}

Our next step is to improve the lower bound in Lemma \ref{B0bounds}. Let
\[
\rho^0_t = N^{-2} A_{W,\psi(t)} - \ep^{7/6}.
\]
On the event
%
%
\begin{equation} \label{Fdef}
F = \bigl\{ \bigl|N^{-2}\mathcal{C}^0_{W,\psi
(s)}\bigr| \ge\rho^0_s \mbox{ for all $s\in
I_{\ep,t}$} \bigr\},
\end{equation}
which has probability
tending to 1 as $\ep\to0$ by Lemma \ref{B0bounds}, $\mathcal
{C}^0_{W,\psi(s)}$ can be coupled with a process $\mathcal{B}^0_{\psi
(s)}$ so that $N^{-2}|\mathcal{B}^0_{\psi(s)}|=\rho^0_s$ and
$\mathcal{C}^0_{W,\psi(s)} \supseteq\mathcal{B}^0_{\psi(s)}$ for
$s\in
I_{\ep,t}$. If for $k \ge1$ $\mathcal{B}^k_{\psi(t)}$ is obtained from
$\mathcal{B}^0_{\psi(t)}$ in the same way as $\mathcal{C}^k_{W,\psi(t)}$
is obtained from $\mathcal{C}^0_{W,\psi(t)}$, then, on $F$, $\mathcal
{C}^k_{W,\psi(s)} \supseteq\mathcal{B}^k_{\psi(s)}$ for $s\in
I_{\ep,t}$. For $k\ge1$ let\looseness=-1
\[
\rho^k_s = N^{-2} \bigl|\mathcal{B}^k_{\psi(s)}\bigr|.
\]\looseness=0
We begin with the case $k=1$. For $f_0(t)=g_0(t)-\ep^{7/6}$, where
$g_0$ is as in (\ref{gdef2}), let
%
%
\begin{equation}\label{f1eq2}
f_1(t) = 1 - \bigl(1-f_0(t)\bigr)
\exp\biggl(-\int_{\log(3\ep)}^t \frac{(t-s)^2}{2}f_0(s)\,ds\biggr).
\end{equation}
\begin{lemma}\label{f1lb}
For any $t<\infty$ there is an $\ep_0=\ep_0(t)>0$ so that for $0 <
\ep< \ep_0$ and any $\delta>0$,
\[
\limsup_{N\to\infty} P\Bigl[\inf_{s\in I_{\ep,t}}
\bigl(N^{-2}C^1_{W,\psi(s)}-f_1(s)\bigr) < -\delta\Bigr] \le P(M<\ep
^{1/3})+\ep^{1/12}.
\]
\end{lemma}
\begin{pf}
As in Lemma \ref{compare1}, if $x\notin\mathcal{B}^0_{\psi(t)}$,
then $x\notin\mathcal{B}^1_{\psi(t)}$ if and only if no generation
1 center is born in the space--time cone
\[
K_{x,t}^\ep\equiv\bigl\{(y,s)\in\Gamma(N)\times[W,\psi(t)]\dvtx
|y-x| \le\bigl(\psi(t)-s\bigr)/\sqrt{2\pi}\bigr\}.\vadjust{\goodbreak}
\]

Conditioning on $\mathcal{G}^0_t=\sigma\{ \mathcal{B}^0_{\psi(s)} \dvtx
s\in
I_{\ep,t}\}$, the locations of generation 1 centers in $\mathcal
{B}^1_{t}$ is a Poisson point process on $\Gamma(N) \times
[W,\psi(t)]$ with intensity
\[
N^{-2} \times|\mathcal{B}^0_{s}|N^{-\alpha} = \rho^0_{\psi
^{-1}(s)} N^{-\alpha}.
\]
Using this and then changing variables $s=\psi(r)$, where
$\psi(r)=R+N^{\alpha/3}r$,
\begin{eqnarray*}
P\bigl( x \notin\mathcal{B}^1_{\psi(t)} | \mathcal{G}^0_t
\bigr) &=& (1-\rho^0_t) \exp\biggl(- \int_W^{\psi(t)}
\frac{(\psi(t)-s)^2}{2}\rho^0_{\psi^{-1}(s)} N^{-\alpha}
\,ds\biggr)\\
&=& (1-\rho^0_t) \exp\biggl(-\int_{\log(3\ep)}^t
\frac{(t-r)^2}{2}\rho^0_r \,dr\biggr).
\end{eqnarray*}

Let $E_{x,t} = \{ x \notin\mathcal{B}^1_{t}\}$.
Since $K_{x,t}^\ep$ and $K_{y,t}^\ep$ are disjoint if
$|x-y|>2(t-\log(3\ep))N^{\alpha/3}/\sqrt{2\pi}$, the events\vspace*{1pt}
$E_{x,t}$ and $E_{y,t}$ are conditionally independent given $\mathcal
{G}^0_t$ if
this holds. Define the random variables $Y_x$, $x\in\Gamma(N)$, so
that $Y_x=1$ if $E_{x,t}$
occurs, and $Y_x=0$ otherwise. From (\ref{cmu1})
%
%
\begin{equation} \label{cmu1}
E( Y_x | \mathcal{G}^0_t)
=(1-\rho^0_t)\exp\biggl(-\int_{\log(3\ep)}^t
\frac{(t-s)^2}{2}\rho^0_s \,ds\biggr).
\end{equation}
Using independence of $Y_x$ and $Y_z$ for
$|x-z|>2(t-\log(3\ep))N^{\alpha/3}/\sqrt{2\pi}$, and the fact that
$\{z\dvtx|x-z| \le2(t-\log(3\ep))N^{\alpha/3}/\sqrt{2\pi
}\}$ has area
$2(t-\log(3\ep))^2 N^{2\alpha/3}$,
%
%
\begin{eqnarray}\label{cvar1}
&&\operatorname{var} \biggl( \int_{x\in\Gamma(N)} Y_x \,dx \big|
\mathcal{G}^0_t\biggr)
\nonumber\\
&&\qquad=\int_{x,z\in\Gamma(N)} [E( Y_xY_z |
\mathcal{G}^0_t) -E( Y_x |
\mathcal{G}^0_t) E( Y_z | \mathcal
{G}^0_t)]
\,dx \,dz
\\
&&\qquad\le N^2\cdot2\bigl(t-\log(3\ep)\bigr)^2 N^{2\alpha/3}.\nonumber
\end{eqnarray}
Using Chebyshev's inequality, we see that
%
%
\begin{eqnarray} \label{cch1}
&&P\biggl(\biggl|\int_{x\in\Gamma(N)}
\bigl(Y_x-E( Y_x | \mathcal{G}^0_t)
\bigr)\, dx \biggr|> \frac{\eta}{2} N^2 \Big|
\mathcal{G}^0_t \biggr) \nonumber\\[-8pt]\\[-8pt]
&&\qquad\le
\frac{4\operatorname{var}(\int_{x\in\Gamma(N)} Y_x \,dx
| \mathcal{G}^0_t)}{\eta^2N^4}.\nonumber
\end{eqnarray}
Combining (\ref{cmu1}), (\ref{cvar1}) and (\ref{cch1}) gives
\begin{eqnarray*}
&&P\biggl( \biggl|(1-\rho^1_t) -
(1-\rho^0_t)\exp\biggl(-\int_{\log(3\ep)}^t \frac{(t-s)^2}{2}\rho^0_s
\,ds\biggr) \biggr|
> \frac{\eta}{2} \Big| \mathcal{G}^0_t \biggr) \\
&&\qquad\le\frac
{8(t-\log(3\ep))^2}{\eta^2N^{2-2\alpha/3}}.
\end{eqnarray*}
The same bound holds for the unconditional probability.
By Lemma \ref{B0bounds} if $\eta>0$ and
\[
F_{0,\eta} \equiv\Bigl\{{\sup_{s\in I_{\ep,t}}} |\rho^0_s-f_0(s)|
\le\eta\Bigr\}\qquad
\mbox{then } \lim_{N\to\infty} P(F_{0,\eta}^c) = 0.
\]
Let $\eta'=\eta[1+(t-\log(3\ep))^3/3!]^{-1}/2$. Using
(\ref{f1eq2}) and the fact that for $x,y\ge0$
%
%
\begin{equation}\label{eineq}
|e^{-x}-e^{-y}| = \biggl| \int_x^y e^{-z} \,dz \biggr| \le|x-y|,
\end{equation}
we see that on the event $F_{0,\eta'}$, we have for any $s\in I_{\ep,t}$
\begin{eqnarray*}
&& \biggl|(1-\rho^0_s) \exp\biggl(-\int_{\log(3\ep)}^s
\frac{(s-r)^2}{2}\rho^0_r \,dr\biggr)-\bigl(1-f_1(s)\bigr)\biggr| \\
&&\qquad \le\bigl|(1-\rho^0_s)-\bigl(1-f_0(s)\bigr)\bigr| + \eta' \int_{\log
(3\ep)}^s \frac{(s-r)^2}{2} \,dr\\
&&\qquad\le\eta' + \eta' \frac{(s-\log(3\ep))^3}{3!} \le\frac{\eta}{2}.
\end{eqnarray*}
So for any $s\in I_{\ep,t}$
\begin{eqnarray*}
&& \lim_{N\to\infty} P\bigl( |\rho^1_s - f_1(s)| >
\eta\bigr)\\
&&\qquad\le\lim_{N\to\infty} P(F_{0,\eta'}^c )\\
&&\qquad\quad{}+ \lim_{N\to\infty}P\biggl(\biggl|(1-\rho^1_s
)-(1-\rho^0_s)
\exp\biggl( -\int_{\log(3\ep)}^s \frac{(s-r)^2}{2}\rho^0_r
\,dr\biggr)\biggr|>\frac{\eta}{2}\biggr) \\
&&\qquad= 0.
\end{eqnarray*}
Since $\eta>0$ is arbitrary, the two quantities being compared are
increasing and continuous, and on the event $F$ defined in
(\ref{Fdef}) $N^{-2} C^1_{W,\psi(s)} \ge\rho^1_s$ for $s\in
I_{\ep,t}$,
\begin{eqnarray*}
& & \limsup_{N\to\infty} P\Bigl[\inf_{s\in
I_{\ep,t}}
\bigl(N^{-2}C^1_{W,\psi(s)}-f_1(s)\bigr) < -\delta\Bigr]\\
&&\qquad \le P(F^c)
+\limsup_{N\to\infty} P\Bigl(\sup_{s\in I_{\ep,t}}
|\rho^1_s-f_1(s)| >\delta\Bigr)
\le P(F^c),
\end{eqnarray*}
and the desired conclusion follows from Lemma \ref{B0bounds}.
\end{pf}

To improve this we will let
%
%
\begin{equation}\label{fiter2}
f_{k+1}(t) = 1 - \bigl(1-f_{k}(t)\bigr)
\exp\biggl(-\int_{\log(3\ep)}^t \frac
{(t-s)^2}{2}\bigl(f_k(s)-f_{k-1}(s)\bigr)\,ds\biggr),\hspace*{-38pt}
\end{equation}
and recall from (\ref{fepinteq}) that as $k\uparrow\infty$,
$f_k(t) \uparrow f_\ep(t)$.
\begin{lemma}\label{fklb}
For any $t<\infty$ there is an $\ep_0=\ep_0(t)>0$ so that for $0 <
\ep< \ep_0$ and any $\delta>0$,
\[
\limsup_{N\to\infty} P\Bigl[\inf_{s\in I_{\ep,t}}
\bigl(N^{-2}C_{\psi(s)}-f_\ep(s)\bigr) < -\delta\Bigr] \le P(M<\ep
^{1/3})+\ep^{1/12}.
\]
\end{lemma}
\begin{pf}
Conditioning on $\mathcal{G}^k_t=\sigma\{ \mathcal{B}^j_{\psi
(s)} \dvtx0\le j\le k, s\in I_{\ep,t}\}$, we have
\[
P\bigl( x \notin\mathcal{B}^{k+1}_{\psi(t)} |
\mathcal{G}^k_t \bigr)
= (1-\rho^k_t) \exp\biggl(-\intt\frac{(t-s)^2}{2}
(\rho^k_s-\rho^{k-1}_s) \,ds\biggr).
\]
Let $F_{k,\eta}=\{\sup_{s\in I_{\ep,t}} |\rho^k_s-f_k(s)|\le\eta\}$,
and $\eta'=\eta[1+2(t-\log(3\ep))^3/3!]^{-1}/2$. Using
(\ref{fiter2}) and $|e^{-x}-e^{-y}| \le|x-y|$
for $x,y\ge0$, we see that on the event $G_{k,\eta'}=F_{k,\eta'}
\cap F_{k-1,\eta'}$, for any $s\in I_{\ep,t}$
\begin{eqnarray*}
&&\biggl|(1-\rho^k_t) \exp\biggl(-\int_{\log(3\ep)}^t
\frac{(t-s)^2}{2}
(\rho^k_s-\rho^{k-1}_s) \,ds
\biggr)-\bigl(1-f_{k+1}(t)\bigr)\biggr|
\\[-2pt]
&&\qquad \le\bigl|(1-\rho_t^k)-\bigl(1-f_k(t)\bigr)\bigr| + 2\eta'
\int_{\log(3\ep)}^t \frac{(t-s)^2}{2} \,ds \\[-2pt]
&&\qquad= \eta'+2\eta' \bigl(t-\log(3\ep)\bigr)^3/3\le\eta/2.
\end{eqnarray*}
Bounding the variance as before we can conclude by induction on $k$
that for any $\eta>0$
%
%
\begin{equation}\label{rhokbd}
\lim_{N\to\infty} P\Bigl( {\sup_{s\in I_{\ep,t}}} |\rho^k_s
- f_k(s)| > \eta\Bigr)= 0.
\end{equation}

Next we bound the difference between $f_k(t)$ and $f_\ep(t)$. Let
$G(t)=t^3/3!$ for $t\ge0$ and $G(t)=0$ for $t<0$. If $*k$
indicates the $k$-fold convolution, then for $k\ge1$, using
arguments similar to the ones in the proof of Lemma \ref{V},
$G^{*k}(t)=t^{3k}/(3k)!$ for $t\ge0$ and $G^{*k}(t)=0$ for $t<0$.
Now if $f*G^{*k}(t) = \intt f(t-r) \,dG^{*k}(r)$, $\tilde f_k(\cdot)=
f_k(\cdot+\log(3\ep))$ and $\tilde f_\ep(\cdot)=f_\ep(\cdot+\log
(3\ep))$,
then changing variables $s\mapsto t-r$ in (\ref{fkinteq}) and
(\ref{fepinteq}), and using the inequality in (\ref{eineq}),
\begin{eqnarray*}
& & \bigl|\tilde f_k\bigl(t-\log(3\ep)\bigr)-\tilde f_\ep\bigl(t-\log
(3\ep)\bigr)\bigr| \\
&&\qquad \le \bigl|{\exp}\bigl(-\tilde f_{k-1}*G\bigl(t-\log(3\ep)\bigr)\bigr)-\exp\bigl(-\tilde
f_\ep*G\bigl(t-\log(3\ep)\bigr)\bigr)\bigr|\\
&&\qquad \le |\tilde f_{k-1}-\tilde
f_\ep|*G\bigl(t-\log(3\ep)\bigr).
\end{eqnarray*}
Iterating the above inequality and using
$|\tilde f_\ep(s)-\tilde f_0(s)|=\tilde f_\ep(s)-\tilde f_0(s)\le1$,
%
%
\begin{eqnarray} \label{fgap}
|f_k(t)-f_\ep(t)| &=& \bigl|\tilde f_k\bigl(t-\log(3\ep)\bigr)-\tilde f_\ep\bigl(t-\log
(3\ep)\bigr)\bigr|
\nonumber\\
& \le & |\tilde f_0-\tilde f_\ep|*G^{*k}\bigl(t-\log(3\ep)\bigr)
\\
& \le & G^{*k}\bigl(t-\log(3\ep)\bigr) = \frac{(t-\log(3\ep))^{3k}}{(3k)!},
\nonumber
\end{eqnarray}
where the last equality comes from (\ref{Fconv}).\vadjust{\goodbreak}

Choose $K=K(\ep,t)$ so that $(t-\log(3\ep))^{3K}/(3K)! <\delta/2$.
Since $C_{\psi(t)} \ge C^k_{W,\psi(t)}$ for any $k\ge0$, and
on the event $F$ defined in (\ref{Fdef}), we have $ C^k_{W,\psi(t)}
\ge|\mathcal{B}^k_{\psi(t)}|$, we have
\[
P\Bigl(\inf_{s\in I_{\ep,t} } \bigl(N^{-2} C_{\psi(s)}-f_\ep
(s)\bigr)<-\delta\Bigr)
\le P(F^c) + P\Bigl({\sup_{s\in I_{\ep,t} }}|\rho
^K_s-f_K(s)| > \delta/2\Bigr).
\]
Using (\ref{rhokbd}) and Lemma \ref{B0bounds} we get the result.
\end{pf}

It is now time to get upper bounds on $C_{\psi(s)}$. Recall $g_0(t)$
defined in (\ref{gdef2}), let
$g_{-1}(t)=0$ and for $k \ge1$ let
%
%
\begin{eqnarray}\label{giter}
g_k(t) &=& 1-\bigl(1-g_{k-1}(t)\bigr)\nonumber\\[-8pt]\\[-8pt]
&&{}\times\exp\biggl(-\int_{\log(3\ep)}^t \frac
{(t-s)^2}{2} \bigl(g_{k-1}(s)-g_{k-2}(s)\bigr) \,ds\biggr).\nonumber
\end{eqnarray}
As in the case of $f_k(t)$, the equations above imply
\[
g_k(t) = 1-\bigl(1-g_{0}(t)\bigr)\exp\biggl(-\int_{\log(3\ep)}^t
\frac{(t-s)^2}{2} g_{k-1}(s) \,ds\biggr),
\]
so we have $g_k(t) \uparrow g_\ep(t)$ as $k\uparrow\infty$, where
$g_\ep(t)$ satisfies
\[
g_\ep(t) = 1-\bigl(1-g_{0}(t)\bigr)\exp\biggl(-\int_{\log(3\ep)}^t
\frac{(t-s)^2}{2} g_\ep(s) \,ds\biggr).
\]
\begin{lemma}\label{glb}
For any $t<\infty$ there exists $\ep_0=\ep_0(t)>0$ such that for $0
< \ep< \ep_0$ and any $\delta>0$,
\[
\limsup_{N\to\infty} P\Bigl[\sup_{s\in I_{\ep,t}}
\bigl(N^{-2}C_{\psi(s)}-g_\ep(s)\bigr) > \delta\Bigr] \le
P(M<\ep^{1/3})+\ep^{2/3}.
\]
\end{lemma}
\begin{pf}
$C^0_{W,\psi(t)} \le A^0_{W,\psi(t)}$. If $\phi^0_t = N^{-2}
A^0_{W,\psi(t)} $ is the fraction\vspace*{1pt} of area covered by generation 0
balloons at time $\psi(t)$, generation 1 centers are born at rate
$N^{2-\alpha}\phi^0_{\psi^{-1}(\cdot)}$.
Let $\phi^1_t$ denotes\vspace*{1pt} the fraction of area covered by centers of
generations $\le1$ at time $\psi(t)$, then
using an argument similar to the one for Lemma \ref{f1lb} gives
\[
\lim_{N\to\infty} P\Bigl( \sup_{s\in I_{\ep,t}} \phi^1_s -
g_1(s) > \eta\Bigr) = 0
\]
for any $\eta>0$. Continuing by induction, let $\phi_t^k$ be the
fraction of area covered by centers of generations $0\le j\le k$.
Since (\ref{giter}) and (\ref{fiter2}) are the same except for
the\vadjust{\goodbreak}
letter they use,
then by an argument identical to the one for Lemma \ref{fklb},
%
%
\begin{equation}\label{eq9}
\lim_{N\to\infty} P\Bigl( {\sup_{s\in
I_{\ep,t}}}
|\phi^k_s - g_k(s)| > \eta\Bigr) = 0
\end{equation}
for any
$\eta>0$. Now using an argument similar to the one for (\ref{fgap})
%
%
\begin{equation}\label{eq8}
{\sup_{s\in I_{\ep,t}}} |g_k(s)-g_\ep(s)| \le\frac
{(t-\log(3\ep))^{3k}}{(3k)!} .
\end{equation}
Next we bound the difference between $C^k_{W,\psi(t)}$ and $C_{\psi
(t)}$. Let
$T=S(\ep^{2/3})$, where $S(\cdot)$ is as in (\ref{S}). Using the
coupling between $\mathcal{C}_t$ and $\mathcal{A}_t$,
\[
C_{\psi(t)}- C^k_{W,\psi(t)} \le A_{\psi(t)}-A_{W,\psi(t)}^k.
\]
Using the fact that $EA_{s+t}-EA_{s,s+t}^k$ is nondecreasing in $s$,
the definitions of $W$ and $T$, Markov's inequality, and Lemma
\ref{Abound}, we have for $T'=T+(t-\log(3\ep))N^{\alpha/3}$,
\begin{eqnarray*}
&&
P\biggl(\sup_{s\in I_{\ep,t}}\bigl(C_{\psi(s)}-\mathcal
{C}^k_{W,\psi(s)}
\bigr) > \frac{\delta N^2}{4} \biggr) \\
&&\qquad \le P(W>T) +
P\biggl( A_{T'}-A_{T,T'} > \frac{\delta N^2}{4} \biggr)\\
&&\qquad \le P(M<\ep^{1/3})+ \frac{4}{\delta N^{2}} E( A_{T'}-A_{T,T'})\\
&&\qquad \le P(M<\ep^{1/3}) + \frac{12\ep^{2/3}}{\delta}
\sum_{j=k+1}^\infty
\frac{(t-\log(3\ep))^j}{j!}.
\end{eqnarray*}
Choose $K=K(\ep,t)$ large enough so that $\sum_{j=K+1}^\infty
(t-\log(3\ep))^j/j! < \delta/12$. If we let
\[
F_K =\Bigl\{\sup_{s\in I_{\ep,t}} \bigl(C_{\psi(s)} -
C^K_{W,\psi(s)}\bigr) \le(\delta/4)N^2\Bigr\},
\]
then
\[
P(F_K^c)\le P(M<\ep^{1/3})+\ep^{2/3}.
\]
By the choice of $K$ and (\ref{eq8}), ${\sup_{s\in I_{\ep,t}}}
|g_K(s)-g_\ep(s)|\le\delta/2$.
Combining the last two inequalities and using the fact that
$N^{-2}C^K_{W,\psi(s)} \le\phi^K_s =\break N^{-2}A^K_{W,\psi(s)}$,
\[
P\Bigl(\sup_{s\in I_{\ep,t}} N^{-2}C_{\psi(s)}-g_\ep(s) > \delta
\Bigr)
\le P(F_K^c)+ P\Bigl({\sup_{s\in I_{\ep,t}}} |\phi
^K_s-g_K(s)| > \delta/4\Bigr).
\]
So using (\ref{eq9}) we have the desired result.
\end{pf}

Our next goal is:
\begin{pf*}{Proof of Lemma \ref{h}}
We prove the result in two steps. To begin we consider a function
$h_\ep(\cdot)$ satisfying
$h_\ep(t) = e^t/3$ for $t < \log(3\ep)$.
%
%
\begin{equation}\label{hep}
h_\ep(t)=1-\exp\biggl(-\int_{-\infty}^{\log(3\ep)}
\frac{(t-s)^2}{2} \frac{e^s}{3} \,ds - \int_{\log(3\ep)}^t \frac
{(t-s)^2}{2} h_\ep(s) \,ds\biggr)\hspace*{-32pt}
\end{equation}
for $t\ge\log(3\ep)$, and prove that $h_\ep(\cdot)$ converges to
some $h(\cdot)$ with the desired properties.
\begin{lemma} \label{hepmono}
For fixed $t$, $h_\ep(t)$ in (\ref{hep}) is monotone decreasing in
$\ep$.
\end{lemma}
\begin{pf}
If we change variables $s = t-u$ and integrate by parts, or remember
the first two moments of the
exponential with mean 1, then
%
%
\begin{eqnarray} \label{id1}
\int_{-\infty}^t (t-s) e^s \,ds &=& \int_0^\infty u e^{t-u} \,du =
e^t,
\nonumber\\[-8pt]\\[-8pt]
\int_{-\infty}^t \frac{(t-s)^2}{2} e^s \,ds &=& \int_0^\infty\frac
{u^2}{2} e^{t-u} \,du
= e^t\int_0^\infty u e^{-u} \,du = e^t.\nonumber
\end{eqnarray}
Using $(t-s)^2/2 = (t-r)^2/2 + (t-r)(r-s) + (r-s)^2/2$ now gives the
following identity
%
%
\begin{equation} \label{id}
\int_{-\infty}^r \frac{(t-s)^2}{2} e^s
\,ds = e^r\biggl[\frac{(t-r)^2}{2}+(t-r)+1\biggr].
\end{equation}
Using (\ref{hep}), the inequality $1-e^{-x}\le x$, (\ref{id1}), and
changing variables $s=t-u$,
\begin{eqnarray*}
h_\ep(t)-\frac13 e^t &\le& \int_{\log(3\ep)}^t \frac
{(t-s)^2}{2}\biggl(h_\ep(s)-\frac13 e ^s\biggr) \,ds \\
&=&\int_0^{t-\log(3\ep)} \biggl(h_\ep(t-u)-\frac13 e^{t-u}
\biggr)\frac{u^2}{2} \,du.
\end{eqnarray*}
Applying Lemma \ref{renewalineq} with $\lambda=1$ and $\beta(\cdot
)\equiv0$ to
$h_\ep(\cdot+\log(3\ep))-\exp(\cdot+\log(3\ep))/3$,
\[
h_\ep(t)-\tfrac13 e^t\le0 \qquad\mbox{for any $t\ge\log(3\ep)$}.
\]
This shows that if $0<\ep<\delta<1$, then $h_\delta(t) \ge h_\ep(t)$
for $t\le\log(3\delta)$. To compare the exponentials for $t >
\log(3\delta)$, we note that
\begin{eqnarray*}
&&\int_{\log(3\ep)}^{\log(3\delta)} \frac{(t-s)^2}{2}
\biggl(h_\ep(s)-\frac13 e^s\biggr) \,ds
+\int_{\log(3\delta)}^t \frac{(t-s)^2}{2}\bigl(h_\ep(s)-h_\delta
(s)\bigr) \,ds\\
&&\qquad \le0+ \int_0^{t-\log(3\delta)} \bigl(h_\ep(t-u)-h_\delta
(t-u)\bigr) \frac{u^2}{2} \,ds.
\end{eqnarray*}
Applying Lemma \ref{renewalineq} with $\lambda=1$ and $\beta(\cdot
)\equiv0$ to
$h_\ep(\cdot+\log(3\delta))-h_\delta(\cdot+\log(3\delta))$, we
see that $h_\ep(t)-h_\delta(t)\le0$
for $t\ge\log(3\delta)$.
\end{pf}
\begin{lemma}
$h(t) =\lim_{\ep\to0} h_\ep(t)$ exists. If $h \not\equiv0$ then
$h$ has properties \textup{(a)--(d)} in Lemma \ref{h}.
\end{lemma}
\begin{pf}
Lemma \ref{hepmono} implies that the limit exists. Since $0\le
h_\ep(t)\le e^t/3$, $ 0\le h(t)\le e^t/3$ and so $\lim_{t\to
-\infty} h(t)=0$. To show that
%
%
\begin{equation}\label{hsatint}
h(t) = 1 - \exp\biggl( -\int_{-\infty}^t \frac{(t-s)^2}{2} h(s) \,ds \biggr),
\end{equation}
we need to show that as $\ep\to0$
%
%
\begin{equation}\label{eq10} \int_{\log(3\ep)}^t\frac{(t-s)^2}{2}
h_\ep(s) \,ds
\to\int_{-\infty}^t \frac{(t-s)^2}{2} h(s) \,ds.
\end{equation}
Given
$\eta>0$, choose $\delta=\delta(\eta)>0$ so that
\[
\delta\bigl[1+\bigl(t-\log(3\delta)\bigr)+\bigl(t-\log(3\delta)\bigr)^2/2\bigr] <
\eta/4.
\]
By bounded convergence theorem, as $\ep\to0$,
\[
\int_{\log(3\delta)}^t \frac{(t-s)^2}{2} h_\ep(s) \,ds \to
\int_{\log(3\delta)}^t \frac{(t-s)^2}{2} h(s) \,ds.
\]
So we can choose $\ep_0=\ep_0(\eta)$ so that the difference between
the two integrals is at most $\eta/2$ for any $\ep<\ep_0$. Therefore
if $\ep<\ep_0$, then
\begin{eqnarray*}
&& \biggl|\int_{\log(3\ep)}^t \frac{(t-s)^2}{2} h_\ep(s) \,ds
-\int_{-\infty}^t \frac{(t-s)^2}{2} h(s) \,ds\biggr|\\
&&\qquad \le\frac{\eta}{2} + 2\int_{-\infty}^{\log(3\delta)} \frac
{(t-s)^2}{2} \frac13 e^s \,ds.
\end{eqnarray*}
Using the identity in (\ref{id}) we conclude that the second term is
\[
\le2\delta\bigl[1+\bigl(t-\log(3\delta)\bigr)+\bigl(t-\log(3\delta)\bigr)^2/2
\bigr]\le\frac{\eta}{2}.
\]
This shows that (\ref{eq10}) holds, and with (\ref{hep}) and
(\ref{id}) proves (\ref{hsatint}).

To prove $\lim_{t\to\infty} h(t)=1$ note that if $h(\cdot)\not
\equiv0$, then there is an $r$ with $h(r)>0$, and so for $t>r$
\[
\int_{-\infty}^t \frac{(t-s)^2}{2} h(s) \,ds \ge h(r)\int_r^t \frac
{(t-s)^2}{2} \,ds = h(r) \frac{(t-r)^3}{3!} \to\infty
\]
as $t\to\infty$. So in view of (\ref{hsatint}), $h(t)\to1$ as
$t\to\infty$, if $h(\cdot)\not\equiv0$.

The last detail is to show if $h(\cdot) \not\equiv0$, then $h(t)
\in(0,1)$ for all $t$. Suppose, if possible, $h(t_0)=0$.
Equation (\ref{hsatint}) implies $\int_{-\infty}^{t_0} h(s)[(t-s)^2/2]
\,ds=0$, and hence $h(s)=0$ for $s\le t_0$. Changing variables
$s\mapsto t-r$, and using (\ref{hsatint}) again with the inequality
$1-e^{-x} \le x$, imply that for any $t>t_0$
\[
h(t)\le\int_{-\infty}^t \frac{(t-s)^2}{2} h(s) \,ds=\int_0^{t-t_0}
h(t-r) \frac{r^2}{2} \,dr.
\]
Applying Lemma \ref{renewalineq} with $\lambda=1$ and
$\beta(\cdot)\equiv0$ to the function $h(\cdot+ t_0)$, we see that
$h(t)\le0$ for any $t>t_0$. But $h(t)\ge0$ for any $t$, and hence
$h \equiv0$, a~contradiction.
\end{pf}

To complete the proof of Lemma \ref{h} it suffices to show that
$|f_\ep(\cdot)-h_\ep(\cdot)|$ and
$|g_\ep(\cdot)-h_\ep(\cdot)|$ converge to 0 as $\ep
\to
0$. To do this, note that if
\[
h_0(t)=1-\exp\biggl(-\int_{-\infty}^{\log(3\ep)}\frac{(t-s)^2}{2}
\frac{e^s}{3} \,ds\biggr),
\]
then
\[
h_\ep(t)=1-\bigl(1-h_0(t)\bigr)\exp\biggl(-\int_{\log(3\ep)}^t \frac
{(t-s)^2}{2} h_\ep(s) \,ds\biggr),
\]
and so using the inequality $ |e^{-x}-e^{-y}|\le|x-y|$ for $x,y\ge
0$,
\[
|h_\ep(t)-g_\ep(t)| \le|h_0(t)-g_0(t)|
+\int_{\log(3\ep)}^t \frac{(t-s)^2}{2}|h_\ep(s)-g_\ep
(s)| \,ds.
\]
Using the inequality $0\le e^{-x}-1+x\le x^2/2$ and the
identity in (\ref{id}),
\begin{eqnarray*}
|h_0(t)-g_0(t)| &\le&\frac12 \biggl[\ep+\ep\bigl(t-\log(3\ep)\bigr)+\ep
\frac{(t-\log(3\ep))^2}{2}\biggr]^2\\
& \le & \frac32 \ep^2\biggl[1+\bigl(t-\log(3\ep)\bigr)^2+\frac{(t-\log(3\ep
))^4}{4}\biggr].
\end{eqnarray*}
Applying Lemma \ref{renewalineq} with $\lambda=1$ and
$\beta(t)=1+t^2+t^4/4$ to the function
\[
\bigl|h_\ep\bigl(\cdot+\log(3\ep)\bigr)-g_\ep\bigl(\cdot+\log(3\ep)\bigr)\bigr|,
\]
we have $|h_\ep(t)-g_\ep(t)| \le(3\ep^2/2)\beta*V(t-\log(3\ep))$,
where $V(\cdot)$ is as in Lem\-ma~\ref{V}. Using $\lambda=1$ in the
expression of $V(\cdot)$ and Lemma \ref{conv},
\begin{eqnarray*}
\beta*V(t) &=&\beta(t)+\intt\beta(t-s) V'(s) \,ds \\
&=& \sum_{k=0}^\infty\biggl[\frac{t^{3k}}{(3k)!}+2\frac{t^{3k+2}}{(3k+2)!}
+6\frac{t^{3k+4}}{(3k+4)!}\biggr] \le6e^t.
\end{eqnarray*}
So $|h_\ep(t)-g_\ep(t)| \le(3\ep^2/2) \cdot6\exp(t-\log(3\ep
))$, and so
\[
{\sup_{s\in I_{\ep,t}}} |h_\ep(s)-g_\ep(s)| \le6\ep e^t/2.
\]
Repeating the argument for $f_\ep(\cdot)$,
and noting that $|h_0(t)-f_0(t)|=|h_0(t)-g_0(t)|+\ep^{7/6}$,
\[
{\sup_{s\in I_{\ep,t}}} |h_\ep(s)-f_\ep(s)| \le
\biggl(6\frac32 \ep^2+\ep^{7/6}\biggr) \exp\bigl(t-\log(3\ep)\bigr)
= \biggl(\frac13 \ep^{1/6}+ 3\ep\biggr) e^t.
\]
This completes the second step and we have proved Lemma \ref{h}.
\end{pf*}

Now we have all the ingredients to prove Theorem \ref{th3}.
\begin{pf*}{Proof of Theorem \ref{th3}}
Let $h(\cdot)$ be as in Lemma \ref{h}. Choose $\ep\in(0,\delta/6)$
small enough so that
\[
{\sup_{s\in I_{\ep,t}}} |g_\ep(s)-h(s)| <\delta/2,\qquad
{\sup_{s\in I_{\ep,t}}} |f_\ep(s)-h(s)| <\delta/2.
\]
Let $D=\{M\le3\ep N^{2-2\alpha/3}\}$. On the event $D$,
$W=\psi(\log(3\ep))>0$. So
%
%
\begin{eqnarray} \label{eq11}
&&
P\Bigl(\sup_{s\le t} \bigl|N^{-2}C_{\psi(s)}-h(s)\bigr|>\delta
\Bigr)\nonumber\\
&&\qquad\le P(D^c) + P\bigl(N^{-2}C_W+h(\log(3\ep))>\delta\bigr)
\nonumber\\[-8pt]\\[-8pt]
&&\qquad\quad{} + P\Bigl(\sup_{s\in I_{\ep,t}}
\bigl(N^{-2}C_{\psi(s)}-h(s)\bigr)>\delta\Bigr) \nonumber\\
&&\qquad\quad{} +
P\Bigl(\inf_{s\in I_{\ep,t}} \bigl(N^{-2}C_{\psi(s)}-h(s)\bigr) <
-\delta\Bigr).\nonumber
\end{eqnarray}
To estimate the second term in (\ref{eq11}) note that $h(\log(3\ep
))\!\le\!(1/3)\exp(\log(3\ep))\!<\delta/2$ and
\[
P(N^{-2}C_W>\delta/2) \le P\bigl(A_W>(\delta
/2)N^2\bigr)\to0
\]
as $N\to\infty$ by Lemma \ref{ALXbd}. To estimate the third term in
(\ref{eq11}) we use Lemma \ref{glb} to get
\begin{eqnarray*}
&& \limsup_{N\to\infty}
P\Bigl( \sup_{s\in I_{\ep,t}} \bigl(N^{-2}C_{\psi(s)}-h(s)\bigr)
>\delta\Bigr)\\
&&\qquad \le\limsup_{N\to\infty} P\Bigl( \sup_{s\in I_{\ep,t}}
\bigl(N^{-2}C_{\psi(s)}-g_\ep(s)\bigr)>\delta/2 \Bigr) \\
&&\qquad\le
P(M<\ep^{1/3})+\ep^{2/3}.
\end{eqnarray*}
For the fourth term in (\ref{eq11}) use Lemma \ref{fklb} to get
\begin{eqnarray*}
&& \limsup_{N\to\infty} P\Bigl(\inf_{s\in I_{\ep,t}}
\bigl(N^{-2}C_{\psi(s)}-h(s)\bigr)<-\delta\Bigr)\\
&&\qquad \le\limsup_{N\to\infty} P\Bigl(\inf_{s\in I_{\ep,t}}
\bigl(N^{-2}C_{\psi(s)}-f_\ep(s)\bigr)<-\delta/2\Bigr) \\
&&\qquad\le
P(M<\ep^{1/3})+\ep^{1/12}.
\end{eqnarray*}
Letting $\ep\to0$, we see that for any $\delta>0$,
%
%
\begin{equation}\label{eq2}
\lim_{N\to\infty} P\Bigl(\sup_{s\in I_{\ep,t}}
\bigl|N^{-2}C_{\psi(s)}-h(s)\bigr|>\delta\Bigr)=0.
\end{equation}

It remains to show that $h(\cdot)\not\equiv0$. Let $\ep, \gamma$
be such that
\[
P[M\le(1+\gamma)\ep^{1/3}] + 11 \frac{\ep^{1/3}}{\gamma} <1.
\]
Fix any $\eta>0$ and let $t_0=\log(3\ep(1+\gamma)+3\eta)$. Using
Lemmas \ref{ALXbd} and \ref{tausigma}
\begin{eqnarray*}
&&\limsup_{N\to\infty} P\bigl(N^{-2} C_{\psi(t_0)} < \ep\bigr) \\
&&\qquad=
\limsup_{N\to\infty} P\bigl(\tau(\ep)>\psi(t_0)\bigr)\\
&&\qquad\le\limsup_{N\to\infty} P\bigl[\tau(\ep) > \sigma\bigl(\ep(1+\gamma)\bigr)\bigr]
+ \limsup_{N\to\infty} P\bigl[\sigma\bigl(\ep(1+\gamma)\bigr)>\psi(t_0)\bigr]\\
&&\qquad\le\limsup_{N\to\infty} P\bigl[\tau(\ep) >
\sigma\bigl(\ep(1+\gamma)\bigr)\bigr]\\
&&\qquad\quad{}
+ \limsup_{N\to\infty} P\bigl(|N^{-2} A_{W_{\ep(1+\gamma
)+\eta}} - \ep(1+\gamma)-\eta| > \eta\bigr)\\
&&\qquad \le P[M\le(1+\gamma)\ep^{1/3}] + 11\frac{\ep^{1/3}}{\gamma} < 1.
\end{eqnarray*}
But if $h(t_0)=0$, we get a contradiction to (\ref{eq2}). This proves
$h(\cdot)\not\equiv0$.
\end{pf*}

\section{Asymptotics for the cover time}\label{sec5}

\mbox{}

\begin{pf*}{Proof of Theorem \ref{th4}}
Theorem \ref{th3} gives a lower bound on the area covered whcih
implies that if $\delta>0$ and $N$ is large, then with high probability
the number of centers in $\mathcal{C}_{\psi(0)}$
dominates a Poisson random variable with mean $\lambda(\delta)
N^{2-(2\alpha/3)}$, where
\[
\lambda(\delta) = \int_{-\infty}^0 \bigl(h(s)-\delta\bigr)^+ \,ds.
\]
If $\delta_0$ is small enough, $\lambda_0\equiv\lambda(\delta_0)
>0$. Dividing the torus into disjoint squares of size $\kappa
N^{\alpha/3} \sqrt{\log N}$, where $\kappa$ is a large constant, the
probability that a given square is vacant is
$\exp(-\lambda_0\kappa^2 \log N)$. If $\kappa\sqrt{\log N} \ge1$,
the number of squares is $\le N^{2-(2\alpha/3)}$. So if
$\lambda_0\kappa^2 \ge2$, then with high probability none of our
squares is vacant. Thus even if no more births of new centers occur
then the entire square will be covered by a time
$\psi(0)+O(N^{\alpha/3} \sqrt{\log N})$.
\end{pf*}

%

%
\printaddresses


\begin{thebibliography}{7}

\bibitem[\protect\citeauthoryear{Aldous}{2007}]{Ald07}
%
\begin{bmisc}[mr]
\bauthor{\bsnm{Aldous},~\bfnm{D.~J.}\binits{D.~J.}}
(\byear{2007}).
\bhowpublished{When knowing early matters: Gossip, percolation, and Nash equilibria.
Preprint. Available at
\href{http://www.stat.berkeley.edu/users/aldous/Unpub/gossip.pdf}{http://www.stat.berkeley.edu/users/aldous/Unpub/}
\href{http://www.stat.berkeley.edu/users/aldous/Unpub/gossip.pdf}{gossip.pdf}.}
\end{bmisc}
%
\endbibitem\vadjust{\goodbreak}

\bibitem[\protect\citeauthoryear{Barbour and Reinert}{2001}]{BR01}
%
\begin{barticle}[auto]
\bauthor{\bsnm{Barbour},~\bfnm{A.~D.}\binits{A.~D.}} \AND
\bauthor{\bsnm{Reinert},~\bfnm{G.}\binits{G.}}
(\byear{2001}).
\btitle{Small worlds}.
\bjournal{Random. Struct. Alg.}
\bvolume{19}
\bpages{54--74}.
\end{barticle}
%
\endbibitem

\bibitem[\protect\citeauthoryear{Cannas, Marco and
Montemurro}{2006}]{CanMarMon06}
%
\begin{barticle}[mr]
\bauthor{\bsnm{Cannas},~\bfnm{Sergio~A.}\binits{S.~A.}},
\bauthor{\bsnm{Marco},~\bfnm{Diana~E.}\binits{D.~E.}} \AND
\bauthor{\bsnm{Montemurro},~\bfnm{Marcelo~A.}\binits{M.~A.}}
(\byear{2006}).
\btitle{Long range dispersal and spatial pattern formation in biological
invasions}.
\bjournal{Math. Biosci.}
\bvolume{203}
\bpages{155--170}.
\bid{doi={10.1016/j.mbs.2006.06.005}, mr={2268332}}
\end{barticle}
%
\endbibitem

\bibitem[\protect\citeauthoryear{Cox and Durrett}{1981}]{CoxDur81}
%
\begin{barticle}[mr]
\bauthor{\bsnm{Cox},~\bfnm{J.~Theodore}\binits{J.~T.}} \AND
\bauthor{\bsnm{Durrett},~\bfnm{Richard}\binits{R.}}
(\byear{1981}).
\btitle{Some limit theorems for percolation processes with necessary and
sufficient conditions}.
\bjournal{Ann. Probab.}
\bvolume{9}
\bpages{583--603}.
\bid{mr={0624685}}
\end{barticle}
%
\endbibitem

\bibitem[\protect\citeauthoryear{Durrett}{2005}]{Dur10}
%
\begin{bbook}[mr]
\bauthor{\bsnm{Durrett},~\bfnm{Rick}\binits{R.}}
(\byear{2005}).
\btitle{Probability: Theory and Examples},
\bedition{3rd} ed.
\bpublisher{Duxbury Press}, \baddress{Belmont, CA}.
\end{bbook}
%
\endbibitem

\bibitem[\protect\citeauthoryear{Durrett}{2007}]{D07}
%
\begin{bbook}[auto]
\bauthor{\bsnm{Durrett},~\bfnm{Rick}\binits{R.}}
(\byear{2007}).
\btitle{Random Graph Dynamics}.
\bpublisher{Cambridge Univ. Press}, \baddress{Cambridge}.
\end{bbook}
%
\endbibitem

\bibitem[\protect\citeauthoryear{Filipe and Maule}{2004}]{FilMau04}
%
\begin{barticle}[mr]
\bauthor{\bsnm{Filipe},~\bfnm{J.~A.~N.}\binits{J.~A.~N.}} \AND
\bauthor{\bsnm{Maule},~\bfnm{M.~M.}\binits{M.~M.}}
(\byear{2004}).
\btitle{Effects of dispersal mechanisms on spatio-temporal development of
epidemics}.
\bjournal{J. Theoret. Biol.}
\bvolume{226}
\bpages{125--141}.
\bid{doi={10.1016/S0022-5193(03)00278-9}, mr={2069297}}
\end{barticle}
%
\endbibitem

\bibitem[\protect\citeauthoryear{Kesten}{1986}]{Kes86}
%
\begin{bincollection}[mr]
\bauthor{\bsnm{Kesten},~\bfnm{Harry}\binits{H.}}
(\byear{1986}).
\btitle{Aspects of first passage percolation}.
In \bbooktitle{\'{E}cole D'\'et\'e de Probabilit\'es de {S}aint-{F}lour,
{XIV}---1984}.
\bseries{Lecture Notes in Math.}
\bvolume{1180}
\bpages{125--264}.
\bpublisher{Springer}, \baddress{Berlin}.
\bid{mr={0876084}}
\end{bincollection}
%
\endbibitem

\end{thebibliography}
\end{document}